\pdfoutput=1
\documentclass[11pt,a4paper]{article}
\IfFileExists{latexml.sty}{\usepackage{latexml}}{\newif\iflatexml\latexmlfalse}
\usepackage[T1]{fontenc}
\usepackage[utf8]{inputenc}
\usepackage{lmodern,mathtools,amssymb,amsthm,microtype}
\usepackage[left=26.3mm,right=26.3mm,top=29mm,bottom=27mm,footskip=11mm]{geometry}
\usepackage{enumitem,booktabs,tabularx,array,titlesec,graphicx,tikz}
\usetikzlibrary{arrows.meta,positioning}
\usepackage[font=small,labelfont=bf]{caption}
\usepackage[unicode,colorlinks=true,linkcolor=blue!45!black,citecolor=blue!45!black,urlcolor=blue!45!black]{hyperref}
\hypersetup{pdftitle={Group algebras that are domains but embed in no division ring},pdfauthor={Nidhal Mghirbi and Seth Douglas},pdfsubject={Group-ring domains and obstructions to division embeddings}}
\titleformat{\section}{\normalfont\normalsize\bfseries}{\thesection.}{.55em}{}
\titleformat{\subsection}{\normalfont\normalsize\itshape}{\thesubsection.}{.5em}{}
\titlespacing*{\section}{0pt}{2.3ex plus .6ex}{1.2ex}
\titlespacing*{\subsection}{0pt}{2ex plus .5ex}{.8ex}
\setlist[enumerate,1]{label=\textup{(\roman*)},leftmargin=2em,itemsep=2pt,topsep=4pt}
\newtheorem{theorem}{Theorem}[section]
\newtheorem{lemma}[theorem]{Lemma}
\newtheorem{proposition}[theorem]{Proposition}
\newtheorem{corollary}[theorem]{Corollary}
\newtheorem{maintheorem}{Theorem}

\numberwithin{equation}{section}
\DeclareMathOperator{\rk}{rank}
\DeclareMathOperator{\diag}{diag}

\newcommand{\F}{\mathbb F}
\newcommand{\Z}{\mathbb Z}
\newcommand{\Q}{\mathbb Q}

\newcommand{\one}{\mathbf 1}

\newcommand{\tabsmall}{\iflatexml\else\small\fi}
\allowdisplaybreaks[2]
\begin{document}
\title{Group algebras that are domains but embed in no division ring}
\author{Nidhal Mghirbi, Seth Douglas\\{\small\itshape finite.works}}
\date{}
\maketitle

\begin{abstract}
We answer Problem 1.5 of the Kourovka Notebook for group algebras over fields of every positive characteristic. For each prime \(p\) we construct a finitely presented, torsion-free, word-hyperbolic group \(\Gamma_p\) whose group algebra is a domain over every field but embeds in no division ring when the field has characteristic \(p\). The group is diffuse, which is what makes its group algebras domains, and it contains an element \(w_p\ne1\) that every homomorphism from \(\Gamma_p\) to the units of a division ring of characteristic \(p\) sends to \(1\). The obstruction is a conflict between two estimates for the ranks of certain matrices over the group algebra: finite-field interpolation bounds their total rank from above, while averaging over permutations of free generators bounds it from below unless certain group elements are identified. A second, retractive construction allows us to combine the obstructions for different primes. In one three-generated diffuse group, a single central element of infinite order dies in every division ring of positive characteristic and in every finite-dimensional linear representation over a field. Lean~4 verifies the nonembedding conclusions for the two prime-specific group algebras.
\end{abstract}

\noindent{\small\textit{Keywords}\quad group algebra, division ring, unique product, diffuse group, hyperbolic group\par
\par}

\section{Introduction}\label{sec:intro}
In 1965, L.~A.~Bokut' recorded the following question in the first Kourovka Notebook \cite[Problem~1.5]{r16}. Can a group ring have no zero divisors and yet admit no embedding in a skew field? For a commutative domain, taking the field of fractions always gives such an embedding. Bokut' later asked explicitly whether group algebras of torsion-free groups embed in division rings \cite[p.~122]{r4}. We answer Problem~1.5 affirmatively for group algebras over fields of each prescribed positive characteristic.

\begin{maintheorem}[Hyperbolic examples]\label{thm:A}
For every prime \(p\) there is a finitely presented, torsion-free, word-hyperbolic diffuse group \(\Gamma_p\) such that
\begin{enumerate}
\item \(K[\Gamma_p]\) is a domain for every field \(K\).
\item For every field \(K\) of characteristic \(p\), there is no injective ring homomorphism from \(K[\Gamma_p]\) to any division ring.
\item There is a single \(w_p\ne1\) in \(\Gamma_p\) such that every homomorphism \(\Gamma_p\to D^\times\) sends \(w_p\) to \(1\), for every division ring \(D\) of characteristic \(p\).
\end{enumerate}
\end{maintheorem}

The last part identifies an element in the kernel of every proposed embedding. The binomial \(w_p-1\) is nonzero in the group algebra because distinct group elements form part of a basis. Every map to a division ring of characteristic \(p\) kills it. The prime here is fixed before the group is chosen. At the end of the paper we construct one group and one central binomial with the same property in every positive characteristic. The integral and characteristic-zero embedding questions remain open, as does Mal'cev's question for right-ordered groups \cite[Problem~1.6]{r16}.

\subsection*{Why matrices enter a question about zero divisors}
The domain assertion has a short explanation. Suppose that \(P\) and \(R\) are the supports of two nonzero elements of \(K[H]\). If some \(g\in PR\) occurs just once as \(xy\), its coefficient in the product is the product of two nonzero coefficients. A group in which this happens for every pair of finite nonempty sets \(P,R\) has unique product. Our groups have the stronger property of diffuseness, which implies unique product \cite[Proposition~1.1]{r6}.

Recent torsion-free examples exhibit zero divisors \cite{r12,r17,r21} or a failure of direct finiteness \cite{r22}, and Jaikin-Zapirain gives hyperbolic examples, with torsion, that fail direct finiteness in every positive characteristic \cite{r13}. Either property excludes an embedding, but a domain is directly finite, since \(ab=1\) gives \((ba-1)b=0\) and hence \(ba=1\). For background on embeddings of group algebras see \cite{r1,r11} and \cite[Section~3.2]{r10}.

An embedding would assign to every matrix over the group algebra the rank of its image over the division ring. We construct a group in which the same matrices would have to satisfy two incompatible estimates. Start with a group \(G_p\) and \(N\) square matrices \(Y_i\) of size \(r\) over \(\F_p[G_p]\). For every homomorphism \(\phi:\F_p[G_p]\to D\) to a division ring,
\begin{equation}\label{eq:intro-budget}
              \sum_i\rk_D\phi(Y_i)<Nr/4 .
\end{equation}
This bound alone leaves open whether the base algebra embeds. We enlarge the group so that unique product survives and sufficiently many distinct elements would force a contradictory lower bound on these same ranks.

Here is how the hyperbolic example does this. After elementary operations, each \(Y_i\) consists of a \(1\)-by-\(1\) identity block and a square array of elements \(e-1\), where the \(e\)'s freely generate a subgroup. If their images are distinct and different from \(1\), the entries of the resulting array are distinct and nonzero. One placement of distinct entries can have small rank. Across all placements, however, their average rank is large. The reason is that when a new row and column are exposed, at most one choice for their common corner can spoil the next pivot. We glue a copy \(G^\pi\) of \(G_p\) to a common free group for every permutation \(\pi\) of the entry positions, and each copy uses the same permutation on all \(N\) arrays. Suppose that a homomorphism \(\phi\) to a division ring keeps the images of the letters in every array distinct and different from \(1\). Averaging over the copies then gives one copy on which the two estimates collide:
\begin{equation}\label{eq:intro-collision}
 \frac{Nr}{2}\ \le\ \sum_i\rk_D\phi(Y_i^\pi)\ <\ \frac{Nr}{4}.
\end{equation}
The left inequality is the average of Lemma~\ref{lem:permutation} over all copies, and the right one is \eqref{eq:intro-budget} restricted to the copy \(G^\pi\). Every representation must therefore kill a letter or identify two letters in one array. A commutator in the central free vertex group turns these finitely many possibilities into the single word \(w_p\).

\subsection*{What the construction has to preserve}
There are two difficulties in realizing this idea. The first is algebraic. OpenAI's construction of a hyperbolic group that is not residually finite supplies marked affine lines, triangular relations, and an interpolation identity \cite[Lemmas~2.2--2.3 and (2.13)--(2.14)]{r20}. Earlier work used these ingredients to obstruct approximate permutation models and Linnell division rings, using fixed-point estimates in the former case and an independence hypothesis in the latter. Splitting the vertices of the triangle complex exhibited the free subgroups freely generated by the normalized group elements \(e\) \cite[Theorems~1.2 and~1.4, Proposition~6.4 and Corollary~6.5]{r9}. We reprove everything we use from these sources in the form needed here: the lines in Lemma~\ref{lem:lines}, the triangle relations in \eqref{eq:triangle}--\eqref{eq:outer}, the interpolation bound in Lemma~\ref{lem:interp} and the free bases in Lemma~\ref{lem:freebasis}. At each point outside a small exceptional set, the incident labels split into two blocks. The matrix sum over each block has rank at most one. Interpolation bounds the total line rank by the sum of these pointwise ranks. To make this work even for small \(p\), we interpolate over a large extension of \(\F_p\) and represent its scalars by matrices over \(\F_p\). This leaves the target division ring unrestricted.

The second difficulty is geometric. Adding relations can destroy the property that made the group algebra a domain. We give the base complex a metric for which a lower bound on lengths of essential loops implies diffuseness. The permutation gluing also requires control of how the free subgroups sit inside the base. Coning off one such subgroup at a time proves the coset-action property that Bowditch's combination theorem requires. Convex lifts of the same subgroups supply the acylindricity used in Dahmani's hyperbolicity theorem. The first argument protects the domain property. The second gives hyperbolicity.

A shorter construction explains how little rank information the algebraic obstruction needs. Adjoin two elements to each free subgroup so that conjugation changes only the diagonal of its array. The difference of the two arrays is invertible as soon as one specified group element survives. This already contradicts the upper bound if it happens on every line. Further conjugations force each array to have full rank in a faithful division-ring image. The enlarged group contains \(\Z^2\), and the construction retains a retraction.

\begin{maintheorem}[Retractive examples]\label{thm:B}
For every prime \(p\), \(G_p\) is a retract of a finitely presented diffuse group \(\Lambda_p\). The algebra \(K[\Lambda_p]\) is a domain for every field \(K\) and embeds in no division ring if \(\operatorname{char}K=p\). There is a nonidentity \(v_p\in[\Lambda_p,\Lambda_p]\) killed by every homomorphism \(\Lambda_p\to D^\times\) for every division ring \(D\) of characteristic \(p\). Writing \(F_i\) for the free subgroups just described and \(U_i\cong\Z^2\) for the subgroup generated by the two elements adjoined at \(F_i\), the kernel is
\[
           \mathop{\ast}_{i,\,gF_i\in G_p/F_i}gU_i g^{-1}.
\]
It is bi-orderable. The group \(\Lambda_p\) is left-orderable if and only if \(G_p\) is left-orderable.
\end{maintheorem}

\begin{figure}[t]
\centering
\iflatexml
\begin{tabular}{@{}c@{}}
\begin{tabular}{cc}
\multicolumn{2}{c}{\fbox{\begin{tabular}{c}\textbf{Base group \(G_p\)} \ (Sections~4--6)\\
  hyperbolic and diffuse, with free subgroups \(F_i\)\\
  every \(\phi\): \(\sum_i\rk_D\phi(Y_i)<Nr\mathord{/}4\) \ (Proposition~\ref{prop:base})\end{tabular}}}\\
\(\swarrow\) \emph{attach shears} & \emph{glue permutation copies} \(\searrow\)\\
\fbox{\begin{tabular}{c}\(\Lambda_p=G_p\ast_{F_i}H_i\), \(H_i=\Z^2\rtimes F_i\) \ (Section~2)\\
  if \(s_i\ne1\) then \(\rk_D\rho(Y_i)\ge1+\tfrac{r-1}{2}\)\\
  \textbf{Theorem~\ref{thm:B}}: retract of \(G_p\)\end{tabular}}
&\fbox{\begin{tabular}{c}\(\Gamma_p\): copies \(G^\pi\) glued to \(A=\mathop{\ast}_iA_i\) \ (Section~3)\\
  some copy has \(\sum_i\rk_D\phi(Y_i^\pi)\ge Nr\mathord{/}2\)\\
  \textbf{Theorem~\ref{thm:A}}: hyperbolic\end{tabular}}\\
\(\downarrow\) \emph{all primes} &\\
\fbox{\begin{tabular}{c}\(E\supset\bigoplus_p\Lambda_p\) \ (Section~8)\\
  \([v_p,x_1]=z\) for every prime \(p\)\\
  \textbf{Theorem~\ref{thm:D}}: \(z\) is central\end{tabular}} &
\end{tabular}
\end{tabular}
\else
\begin{tikzpicture}[font=\small,
  box/.style={draw, rounded corners=2pt, align=center, inner sep=4pt},
  arr/.style={-{Stealth[length=2.2mm]}, thick}]
\node[box, minimum width=84mm] (G) at (0,0) {\textbf{Base group \(G_p\)} \ (Sections~4--6)\\
  hyperbolic and diffuse, with free subgroups \(F_i\)\\
  every \(\phi\): \(\sum_i\rk_D\phi(Y_i)<Nr/4\) \ (Proposition~\ref{prop:base})};
\node[box, minimum width=68mm] (L) at (-4.15,-2.75) {\textbf{\(\Lambda_p=G_p\ast_{F_i}H_i\)}, \(H_i=\Z^2\rtimes F_i\) \ (Section~2)\\
  if \(s_i\ne1\) then \(\rk_D\rho(Y_i)\ge1+\tfrac{r-1}{2}\)\\
  \textbf{Theorem~\ref{thm:B}}: retract of \(G_p\)};
\node[box, minimum width=68mm] (P) at (4.15,-2.75) {\textbf{\(\Gamma_p\)}: copies \(G^\pi\) glued to \(A=\mathop{\ast}_iA_i\) \ (Section~3)\\
  some copy has \(\sum_i\rk_D\phi(Y_i^\pi)\ge Nr/2\)\\
  \textbf{Theorem~\ref{thm:A}}: hyperbolic};
\node[box, minimum width=68mm] (E) at (-4.15,-5.35) {\textbf{\(E\supset\bigoplus_p\Lambda_p\)} \ (Section~8)\\
  \([v_p,x_1]=z\) for every prime \(p\)\\
  \textbf{Theorem~\ref{thm:D}}: \(z\) is central};
\draw[arr] ([xshift=-16mm]G.south) -- node[above left, font=\footnotesize]{attach shears} (L.north);
\draw[arr] ([xshift=16mm]G.south) -- node[above right, font=\footnotesize]{glue permutation copies} (P.north);
\draw[arr] (L.south) -- node[right, font=\footnotesize]{all primes} (E.north);
\end{tikzpicture}
\fi
\caption{How the groups are built. Every arrow preserves diffuseness, and hence the domain property (Lemma~\ref{lem:extension} and Bowditch's theorems); the right-hand arrow also preserves hyperbolicity (Bowditch, Dahmani). In each of \(\Lambda_p\) and \(\Gamma_p\) the lower bound in the box contradicts the upper bound in \(G_p\) unless a specified letter dies.}
\label{fig:roadmap}
\end{figure}
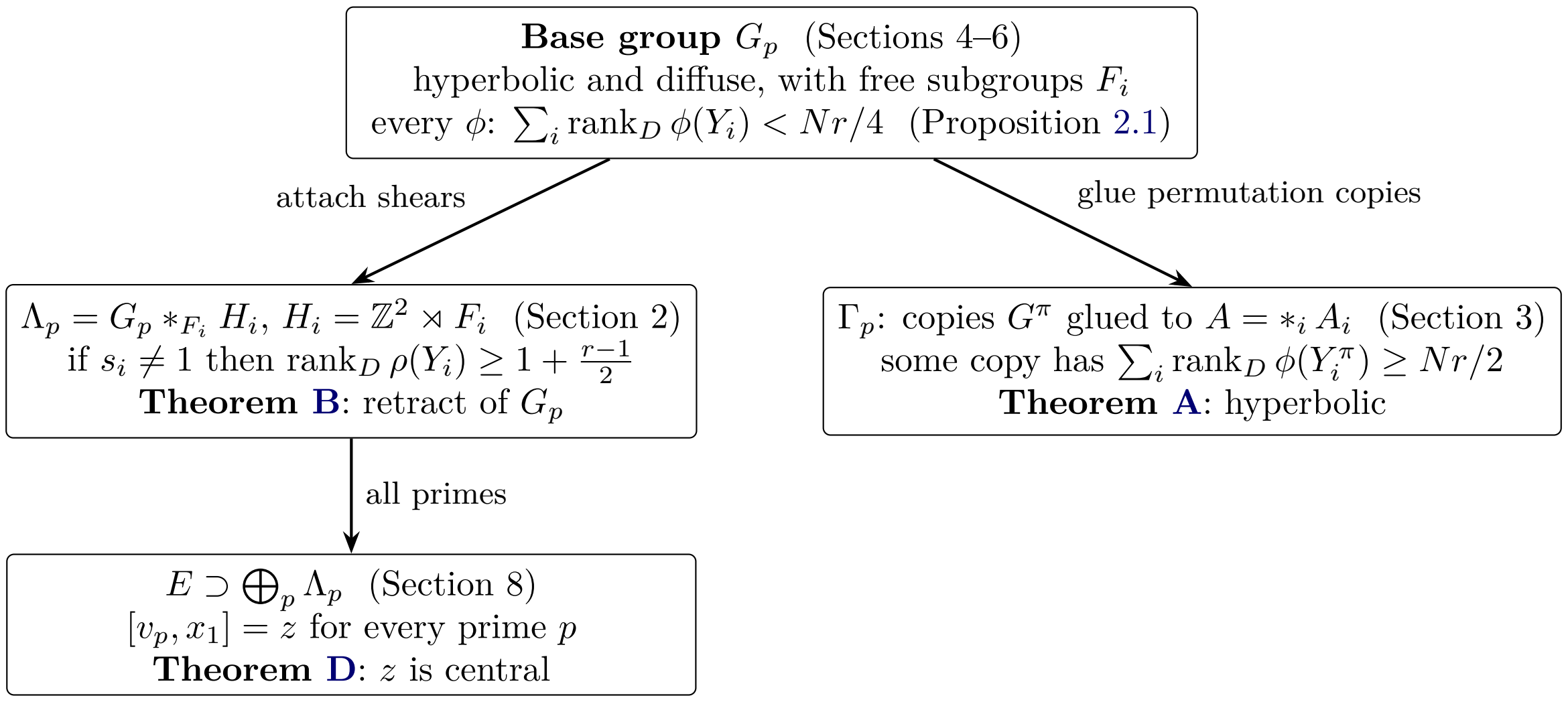

We first deduce both theorems from precise properties of the base (Proposition~\ref{prop:base}). The shear argument also gives a general attachment theorem over any field (Theorem~\ref{thm:C}). Sections~4--6 construct the base and prove its properties. Section~7 gives finite presentations and explains what is known about embeddings of the base itself. In Section~\ref{sec:uniform} the words \(v_p\) become commutators with one common value, a central element of infinite order. Appendix~A identifies the two nonembedding conclusions checked in Lean and their relation to the presentations used here. Figure~\ref{fig:roadmap} shows how the groups are built from one another, and Table~\ref{tab:notation} collects the main notation.

\begin{table}[t]
\centering\tabsmall
\renewcommand{\arraystretch}{1.12}
\begin{tabularx}{\linewidth}{@{}lX@{}}
\toprule
\(G_p\), \(X\) & the base group and its triangle complex, \(G_p=\pi_1X\) (Sections~4--6)\\
\(N\), \(r\) & the number and the size of the line matrices\\
\(Y_i\), \(M_i\), \(e_{i,uv}\) & the line matrices, their normalized \((r-1)\times(r-1)\) arrays, and the free generators in the arrays \eqref{eq:normal}\\
\(F_i\) & the free subgroup generated by the \(e_{i,uv}\)\\
\(U_i=\langle s_i,\gamma_i\rangle\), \(H_i\) & the shear group \(\Z^2\) and \(H_i=U_i\rtimes F_i\) \eqref{eq:semidirect}\\
\(\Lambda_p\), \(v_p\) & the retractive group and its common word (Theorem~\ref{thm:B})\\
\(A=\mathop{\ast}_iA_i\), \(\Pi\), \(G^\pi\) & the central free group, the permutations of entry positions, and the copies of \(G_p\)\\
\(\Gamma_p\), \(w_p\) & the hyperbolic group and its common word (Theorem~\ref{thm:A})\\
\(E\), \(z\) & the three-generated group and its central element (Theorem~\ref{thm:D})\\
\(D\), \(\rk_D\) & a division ring and rank over it, the dimension of the right column space\\
\bottomrule
\end{tabularx}
\caption{Main notation.}
\label{tab:notation}
\end{table}

Throughout, rings and their homomorphisms are unital. Rank over a division ring means the dimension of the right column space. We write \([a,b]=aba^{-1}b^{-1}\) and \(D^\times\) for the multiplicative group of \(D\).

\section{Forcing rank by conjugation}

We state here what the later construction gives. This permits us to prove the nonembedding assertion before making any choices of lengths, girths or field sizes. A finite set in a group action has an extremal point \(x\) if \(hx,h^{-1}x\) can both belong to the set only when \(hx=x\). An action is diffuse if every finite subset of at least two points has two extremal points. For a group itself we use its action by left multiplication. Bowditch calls a subgroup \(F\) of a torsion-free group \(G\) concentrated if the action of \(G\) on \(G/F\) is diffuse \cite{r6}.

\begin{proposition}[The base group]\label{prop:base}
For each prime \(p\) there are a finitely presented, torsion-free, word-hyperbolic diffuse group \(G=G_p\), integers \(N\ge2,r\ge3\), and \(r\times r\) matrices \(Y_1,\ldots,Y_N\) over \(\F_p[G]\) with the following properties.
\begin{enumerate}
\item Invertible row and column operations put \(Y_i\) in the form
\begin{equation}\label{eq:normal}
\begin{pmatrix}1&\one^{T}\\ \one&(e_{i,uv})\end{pmatrix}
  \ \sim\ \diag(1,M_i),\qquad M_i=(e_{i,uv}-1)_{1\le u,v<r}.
\end{equation}
For each \(i\) the \(e_{i,uv}\) form a free basis of a subgroup \(F_i\).
\item Every homomorphism \(\phi:\F_p[G]\to D\) to a division ring satisfies \(\sum_i\rk_D\phi(Y_i)<Nr/4\).
\item The \(F_i\) are quasiconvex and concentrated in \(G\). They are jointly malnormal, meaning that \(F_i\cap gF_jg^{-1}\ne1\) implies \(i=j\) and \(g\in F_i\).
\end{enumerate}
\end{proposition}

The base and its rank bound are constructed in Section~4. Section~5 proves freeness, diffuseness and joint malnormality. The proof of concentration occupies Section~6. Only the hyperbolic gluing needs part~(iii).

\subsection{Changing the diagonal}
Let \(F\) be free on \(k^2\) elements \(e_{uv}\). We seek relations that control the matrix \(M=(e_{uv}-1)\) in a division-ring image. Give \(e_{uv}\) weight \(d(e_{uv})=\delta_{uv}\), so that \(d:F\to\Z\) is a homomorphism. Let \(U=\langle s,\gamma\mid [s,\gamma]=1\rangle\cong\Z^2\). The group \(F\) acts on \(U\) by fixing \(s\) and sending \(\gamma\) to \(s^{-d(h)}\gamma\) under \(h\in F\). On exponent pairs this is the shear \((a,b)\mapsto(a-d(h)b,b)\). Form
\begin{equation}\label{eq:semidirect}
                         H=U\rtimes F .
\end{equation}
The projection \(H\to F\) is a retraction, \(s\) is central and has infinite order, and the defining conjugation takes the form
\begin{equation}\label{eq:shear}
                    \gamma h\gamma^{-1}=s^{d(h)}h\qquad(h\in F).
\end{equation}
For example, in size two, writing \(e_{11}=a,e_{12}=b,e_{21}=c,e_{22}=d_0\), one conjugation gives
\[
\begin{pmatrix}sa-1&b-1\\c-1&sd_0-1\end{pmatrix}
-\begin{pmatrix}a-1&b-1\\c-1&d_0-1\end{pmatrix}
=(s-1)\begin{pmatrix}a&0\\0&d_0\end{pmatrix}.
\]
In size \(k\) the same subtraction gives \((s-1)\diag(e_{11},\ldots,e_{kk})\). Under any homomorphism \(\rho:H\to D^\times\) with \(\rho(s)\ne1\), its image is invertible. The two matrices on the left have the same rank, since one is conjugate to the other. Writing \(\rho(M)\) for \((\rho(e_{uv})-1)_{uv}\), rank subadditivity gives
\begin{equation}\label{eq:half}
                            k\le 2\rk_D\rho(M).
\end{equation}
This one subtraction will give a word killed by every representation in Theorem~\ref{thm:B}. Iterating the conjugation gives a statement useful even when the upper bound falls short of full rank by only one.

\begin{lemma}[Full rank under a shear]\label{lem:shear}
For every field \(K\) and every injective homomorphism \(K[H]\to D\) to a division ring, \(M=(e_{uv}-1)_{u,v=1}^k\) is invertible over \(D\).
\end{lemma}
\begin{proof}
Work in the division subring \(E\) generated by the image of \(K[H]\). The element \(s\) commutes with the generators of \(E\), hence also with their inverses and with all of \(E\). Its powers are distinct. Write \(\Delta=\diag(e_{11},\ldots,e_{kk})\) and \(C=\Delta^{-1}(M-\Delta)\). Conjugation fixes \(M-\Delta\), whose diagonal is \(-1\), and multiplies \(\Delta\) by \(s\). Thus
\begin{equation}\label{eq:iterated-shear}
           (\gamma^j I)M(\gamma^{-j} I)=\Delta(s^jI+C)\qquad(j\ge0).
\end{equation}
Singularity over \(E\) would give \(C\) a central right eigenvalue \(-s^j\) for every \(0\le j\le k\). Eigenvectors for distinct central eigenvalues are right-linearly independent. Indeed, in a shortest dependence, apply \(C\) and then subtract the dependence multiplied on the right by any one of its eigenvalues. Centrality cancels that term and leaves every other coefficient nonzero. This contradicts minimality. In dimension \(k\) there can be no such \(k+1\) vectors. The inverse over \(E\) also serves over \(D\).
\end{proof}

\subsection{Attaching the shears}
\begin{maintheorem}[A general attachment theorem]\label{thm:C}
Let \(K\) be a field and \(G\) a unique-product group. Suppose that finitely many matrices \(Y_i\) over \(K[G]\), of orders \(r_i\ge2\), can be put in the first form in \eqref{eq:normal} by invertible row and column operations over \(K[G]\). For each \(i\) suppose the interior entries \(e_{i,uv}\) freely generate a subgroup \(F_i\le G\). If every embedding \(\phi:K[G]\hookrightarrow D\) in a division ring satisfies \(\sum_i\rk_D\phi(Y_i)<\sum_i r_i\), then \(G\) is a retract of a unique-product group \(\Gamma\) for which \(K[\Gamma]\) embeds in no division ring. If \(G\) is diffuse, so is \(\Gamma\). If \(G\) is finitely presented, so is \(\Gamma\).
\end{maintheorem}

For each \(F_i\), form the preceding semidirect product \(H_i=U_i\rtimes F_i\), with \(U_i=\langle s_i,\gamma_i\rangle\cong\Z^2\). Take the star of groups with central vertex \(G\), leaves \(H_i\), and edge groups \(F_i\) included at both ends. Its fundamental group is
\begin{equation}\label{eq:star}
                        \Gamma=G\mathop{\ast}_{F_i,\,i\in I}H_i .
\end{equation}
Normal form embeds all the vertex groups \cite[I.5.2, Theorem~11 and Corollary~3]{r23}. The projections of the leaves onto \(F_i\) agree with the identity of \(G\) on the edges and give \(\pi:\Gamma\to G\).

\begin{lemma}\label{lem:free-kernel}
The kernel of \(\pi\) is the free product
\begin{equation}\label{eq:free-kernel}
                  \mathop{\ast}_{i\in I,\,gF_i\in G/F_i}gU_i g^{-1}.
\end{equation}
\end{lemma}
\begin{proof}
Let \(L=\ker\pi\) act on the Bass--Serre tree. The stabilizers of edges and of vertices over \(G\) are trivial. Over a leaf the stabilizers are conjugates of \(U_i\). For a vertex or edge group \(B\), its orbits are \(L\backslash\Gamma/B=G/\pi(B)\) under the natural identification. Consequently the quotient graph has one central vertex, one leaf for each pair \(i,gF_i\), and one edge to that leaf. Taking the fundamental group of this tree of groups gives \eqref{eq:free-kernel} by \cite[I.5.4, Theorem~13]{r23}.
\end{proof}

\subsection{Why the algebra is still a domain}
Each \(U_i\) is diffuse, as the two ends of a finite set in a lexicographic order on \(\Z^2\) are extremal. Free products of diffuse groups are diffuse \cite[Theorem~1.2(2)]{r6}, also for infinitely many factors since a finite set uses only finitely many. Thus \(L\) in \eqref{eq:free-kernel} is diffuse and has unique product.

\begin{lemma}\label{lem:extension}
Both unique product and diffuseness are preserved under extensions.
\end{lemma}
\begin{proof}
Consider \(1\to L\to H\xrightarrow{\pi}Q\to1\). For unique product, given finite nonempty \(P,R\subset H\), choose a uniquely represented product \(ab\) in \(\pi(P)\pi(R)\). Only the fibers \(P_a,R_b\) can contribute to products above \(ab\). Write them as \(xP'\) and \(R'y\), respectively, with \(P',R'\subset L\). A uniquely represented product in \(P'R'\), translated by \(x,y\), has a unique representation in \(PR\).

For diffuseness take a finite \(S\subset H\) with at least two elements. If \(\pi(S)\) has two or more elements, choose two extremal ones in \(\pi(S)\) and over each choose an extremal point of its fiber under \(L\) (the point itself for a singleton). An element which moves both ways inside \(S\) at such a point first has to lie in \(L\), by extremality in \(Q\), and then has to fix that point, by extremality in the fiber. If \(S\) lies in one fiber, identify the fiber with an \(L\)-orbit and choose its two extremal points there. This is Bowditch's extension argument \cite[Theorem~1.2(1)]{r6}.
\end{proof}

\begin{proof}[Proof of Theorem~\ref{thm:C}]
Lemma~\ref{lem:extension} gives unique product for \(\Gamma\) and diffuseness when \(G\) is diffuse. In a faithful division-ring image each \(K[H_i]\) embeds, so Lemma~\ref{lem:shear} makes \(M_i\) invertible. Subtracting the first row of the normalized \(Y_i\) from every other row and then clearing the first row gives \(\diag(1,M_i)\). All \(Y_i\) therefore have full rank, contradicting the hypothesis on the restriction to \(K[G]\). A finite presentation of \(G\) extends to one of \(\Gamma\) by adding two generators per line and the relations in \eqref{eq:line-shear-presentation}.
\end{proof}

\begin{proof}[Proof of Theorem~\ref{thm:B} from Proposition~\ref{prop:base}]
Use \eqref{eq:star} with \(G=G_p\), and call the resulting group \(\Lambda_p\). The kernel calculation and Lemma~\ref{lem:extension} give diffuseness, and the finite presentation is the one used for Theorem~\ref{thm:C}. Thus \(K[\Lambda_p]\) is a domain for every \(K\).

Let \(\rho:\Lambda_p\to D^\times\) be any group homomorphism in characteristic \(p\). The central prime field of \(D\) extends \(\rho\) linearly to \(\F_p[\Lambda_p]\). If every \(s_i\) survived, \eqref{eq:normal} and \eqref{eq:half} would give \(\rk_D\rho(Y_i)\ge1+(r-1)/2\) for every \(i\). Summing contradicts Proposition~\ref{prop:base}(ii) on \(\F_p[G_p]\). At least one \(s_i\) is therefore killed, although the line \(i\) may depend on \(\rho\).

Here we can remove that dependence. The elements \(s_1,\ldots,s_N\) belong to distinct free factors in \eqref{eq:free-kernel}. Set \(u_1=s_1\) and \(u_j=[u_{j-1},s_j]\) for \(2\le j\le N\). If \(u_{j-1}\ne1\), the four syllables in this commutator alternate between the subgroup generated by the earlier factors and the new factor. Thus \(u_j\ne1\) by free-product normal form. On the other hand, killing any \(s_i\) makes \(u_N\) the identity. Take \(v_p=u_N\). It belongs to the commutator subgroup because \(N\ge2\), and it has infinite order because \(\Lambda_p\) is diffuse. Every map to a division ring kills \(v_p-1\) over a field of characteristic \(p\), by restriction to the group basis.

The kernel is bi-orderable by Vinogradov's free-product theorem \cite{r24}. Left orders restrict to subgroups and extend across group extensions. For the latter assertion, declare an element positive if its image in the quotient is positive, or if its image is trivial and it is positive in the kernel. This positive cone is closed under multiplication and selects one of \(h,h^{-1}\) for each \(h\ne1\). Thus \(\Lambda_p\) admits a left order exactly when \(G_p\) does.
\end{proof}

\section{A hyperbolic group from permutations}\label{sec:hyperbolic}

The group \(\Lambda_p\) contains the subgroup
\(\langle s_i,\gamma_i\rangle\cong\mathbb Z^2\) and hence cannot be
word-hyperbolic. Another way to obtain a lower bound on rank uses only
free edge groups. Its starting point is an elementary observation
about ordinary elimination that also works over a division ring.

\subsection{Averaging the rank of an array}

Consider a \(2\times2\) array with first row \(a,b\), second row
\(c,z\), and \(a\ne0\). With \(a,b,c\) fixed there is exactly one choice
of \(z\) for which the second row fails to raise the rank, namely
\(z=ca^{-1}b\). The following proof uses this same observation even
when the square already exposed is singular.

\begin{lemma}[Permutation average]\label{lem:permutation}
Take \(n^2\) pairwise distinct nonzero elements of a division ring \(D\),
where \(n\ge2\). Let \(B^\pi\) range over their placements in an
\(n\times n\) array. Then
\begin{equation}\label{eq:perm-avg}
 \frac{1}{(n^2)!}\sum_\pi \operatorname{rank}_D B^\pi
 \ge (n-1)\left(1-\frac1{2n}\right)\ge\frac{n-1}{2}.
\end{equation}
\end{lemma}

\begin{proof}
Suppose the exposed leading square has rank \(s\). Choose \(s\)
pivot rows and columns giving an invertible minor \(P\). Reveal the
next bordering row and column except their common corner \(z\). On the
pivot rows and columns the enlarged matrix contains
\(\left(\begin{smallmatrix}P&u\\v&z\end{smallmatrix}\right)\).
If \(z\ne vP^{-1}u\), this minor has rank \(s+1\), by elimination.
For \(s=0\), the same statement reads \(z\ne0\). Thus at most one
choice of the new corner can fail to increase the rank of the leading
square, regardless of its earlier pivots.

For a uniform permutation, at stage \(k\) there are exactly
\(n^2-k^2+1\) possible values for the corner of the leading
\(k\times k\) square, conditional on the exposed entries. For
\(1\le k\le n-1\) this count is at least \(2n\). Rank never decreases
on bordering, so the expected increase at each of these stages is
at least \(1-1/(2n)\). Adding the expectations proves the first bound,
even without using the last stage. The second bound follows at once.
\end{proof}

\subsection{Making all placements appear in one group}

Put \(n=r-1\), \(U=\{1,\ldots,n\}^2\), and
\(\Pi=\operatorname{Sym}(U)\). Let \(A=\mathop{*}_{i=1}^N A_i\), where
\(A_i\) is free on the \(n^2\) letters \(a_{i,uv}\). Use \(A\) as one
vertex group and take one marked copy \(G^\pi\) of \(G_p\) for each
\(\pi\in\Pi\). For each pair \((\pi,i)\) join \(G^\pi\) to \(A\) by an
edge with group \(F_i\) and maps
\begin{equation}\label{eq:perm-edge-maps}
 e_{i,uv}\longmapsto e_{i,uv}^{\pi}\in G^\pi,\qquad
 e_{i,uv}\longmapsto a_{i,\pi(u,v)}\in A_i .
\end{equation}
A permutation here may exchange arbitrary entry positions. There are
\(N\) parallel edges from \(A\) to each \(G^\pi\).

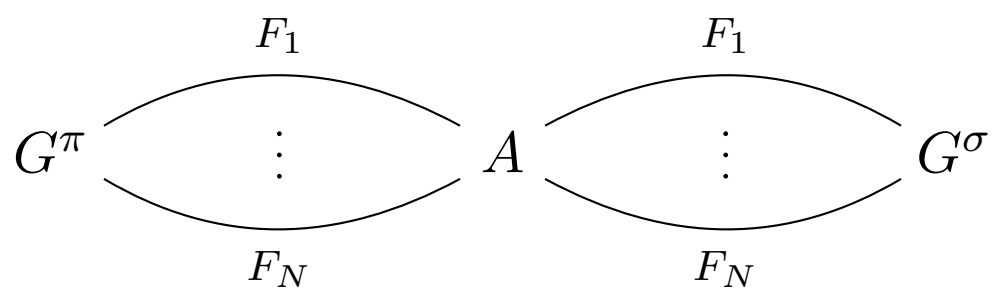
\begin{figure}[htbp]
 \centering
 \begin{tikzpicture}[x=1pt,y=1pt]
  \node at (12,25) {\(G^\pi\)};
  \node at (97,25) {\(A\)};
  \node at (182,25) {\(G^\sigma\)};
  \draw (22,30) .. controls (43.33,42.67) and (65.67,42.67) .. (89,30);
  \draw (22,20) .. controls (43.33,7.33) and (65.67,7.33) .. (89,20);
  \draw (105,30) .. controls (128.33,42.67) and (150.67,42.67) .. (172,30);
  \draw (105,20) .. controls (128.33,7.33) and (150.67,7.33) .. (172,20);
  \node at (55,47) {\(\scriptstyle F_1\)};
  \node at (55,3) {\(\scriptstyle F_N\)};
  \node at (139,47) {\(\scriptstyle F_1\)};
  \node at (139,3) {\(\scriptstyle F_N\)};
  \foreach \y in {20.5,24.5,28.5} {\fill (55,\y) circle (.45pt); \fill (139,\y) circle (.45pt);}
 \end{tikzpicture}
 \caption{Two of the copies attached to the free group \(A\).
 At the central vertex the two edges labelled \(F_i\) have the same
 image \(A_i\), with their bases permuted by \(\pi\) and \(\sigma\).}
 \label{fig:permutation-graph}
\end{figure}

Let \(\Gamma_p\) be the fundamental group of this graph of groups.
Choose one label \(i_*\). The edges \((\pi,i_*)\) form a spanning tree.
In the associated presentation,
\begin{equation}\label{eq:perm-stable-letters}
 t_{\pi,i}e^\pi_{i,uv}t_{\pi,i}^{-1}=a_{i,\pi(u,v)},\qquad
 t_{\pi,i_*}=1.
\end{equation}
A single conjugator serves the whole array for each \((\pi,i)\).
The edge maps are injective on the displayed free bases, so
Bass--Serre normal form embeds every vertex group
\cite[I.5.2, Theorem 11 and Corollary 3, and I.5.3, Theorem 12]{r23}.
The graph is finite, the vertex groups finitely presented and the
edge groups finitely generated. This gives a finite presentation
with \((N-1)(n^2)!\) stable letters. A finite-order element fixes
a vertex of the Bass--Serre tree, since the action has no inversions.
The vertex groups are torsion-free, so \(\Gamma_p\) is torsion-free
\cite[I.4.3, Proposition 19]{r23}.

We show that every homomorphism
\(\rho:\Gamma_p\to D^\times\), with \(\operatorname{char}D=p\), kills
at least one word in a fixed list in the central vertex. Order \(U\)
lexicographically and write \(a_{i,\alpha}\) for its entry letters.
List the letters \(a_{i,\alpha}\) in the order of \((i,\alpha)\),
followed by the words
\(a_{i,\alpha}a_{i,\beta}^{-1}\) for \(\alpha<\beta\) in the order of
\((i,\alpha,\beta)\). Write the resulting list as
\(\ell_1,\ldots,\ell_M\). All these words are nontrivial in \(A\).

Suppose that \(\rho(\ell_j)\ne1\) for every \(j\), and let
\(\phi:\mathbb F_p[\Gamma_p]\to D\) be the unital homomorphism
obtained by extending \(\rho\) linearly over the central prime
field. For each fixed \(i\), the differences
\(b_{i,uv}=\phi(a_{i,uv})-1\) are distinct and nonzero. Write
\(B_i^\pi=(b_{i,\pi(u,v)})\). Equation~\eqref{eq:perm-stable-letters}
gives
\[
 \phi(e^\pi_{i,uv})-1
 =\rho(t_{\pi,i})^{-1} b_{i,\pi(u,v)}\rho(t_{\pi,i}).
\]
On matrices this is multiplication on the two sides by inverse
scalar diagonal matrices, hence it preserves rank. Consequently
(2.1) gives
\begin{equation}\label{eq:perm-rank-match}
 \operatorname{rank}_D\phi(Y_i^\pi)=1+\operatorname{rank}_D B_i^\pi .
\end{equation}
Apply Lemma~\ref{lem:permutation} to each line and sum over \(i\).
Linearity of the average gives a lower bound for the total on a
single copy, even though the same permutation acts on every line.
\begin{equation}\label{eq:perm-lower}
 \max_\pi\sum_i\operatorname{rank}_D\phi(Y_i^\pi)
 \ge \frac{1}{|\Pi|}\sum_\pi\sum_i\operatorname{rank}_D\phi(Y_i^\pi)
 \ge N\left(1+\frac{n-1}{2}\right)=\frac{Nr}{2}.
\end{equation}
On the other hand, restriction to every \(G^\pi\) gives, by
Proposition~2.1(ii),
\begin{equation}\label{eq:perm-upper}
 \sum_i\operatorname{rank}_D\phi(Y_i^\pi)<Nr/4 .
\end{equation}
The proposition applies to every unital homomorphism. This
contradiction shows that at least one
\(\ell_j\) lies in \(\ker\rho\).

We combine the list into a single nontrivial word. For any nontrivial
\(x,y\in A\), some \(h\in A\) satisfies
\([x,hyh^{-1}]\ne1\), where \([x,y]=xyx^{-1}y^{-1}\). If \(x\) and
\(y\) do not commute, take \(h=1\). If they commute, their axes in
the Cayley tree of \(A\) coincide. The stabilizer of that line
acts faithfully on it by translations, since the action on the
Cayley tree is free and without inversions. The stabilizer is cyclic,
whereas \(A\) is nonabelian, so some free basis letter \(h\) does not
preserve the line. The axis of \(hyh^{-1}\) then differs from that
of \(x\), and the elements do not commute.

Set \(w_1=\ell_1\). Given \(w_j\ne1\), choose \(h_j\) as the
first element \(h\) in the list of \(1\) followed by the ordered
basis letters of \(A\) with \([w_j,h\ell_{j+1}h^{-1}]\ne1\), and set
\[
 w_{j+1}=[w_j,h_j\ell_{j+1}h_j^{-1}].
\]
Thus \(w_p=w_M\) is nontrivial in \(A\) and hence in \(\Gamma_p\).
If \(\rho(\ell_1)=1\), then \(\rho(w_1)=1\). If
\(\rho(\ell_j)=1\) for \(j>1\), the second entry in the commutator
defining \(w_j\) maps to \(1\), so \(\rho(w_j)=1\). In either case
all later \(w_k\) map to \(1\). Thus \(\rho(w_p)=1\) for every
\(\rho\) under consideration. A division embedding of
\(K[\Gamma_p]\) with \(\operatorname{char}K=p\) would give a faithful
such homomorphism on the group elements. This proves
Theorem~\ref{thm:A}(ii) and (iii).

\subsection{The two geometric conclusions}

\begin{proof}[Completion of Theorem~\ref{thm:A} from Proposition~2.1]
First consider diffuseness. A free factor of a free group is
concentrated. Indeed, the free factors are diffuse, and in their
free product the trivial edge group is concentrated at both ends.
Bowditch's graph-of-groups theorem, including its assertion for a
concentrated subgroup of a vertex group, applies to each factor,
which is concentrated in itself
\cite[Theorems 1.2(3) and 1.3(2)]{r6}. This applies equally in
infinite rank. In the permutation graph each edge group is free and
diffuse. Its image is concentrated in \(A\) by the free-factor
observation and in \(G^\pi\) by Proposition~2.1(iii). Bowditch's
theorem now gives diffuseness of the torsion-free group \(\Gamma_p\).
Its hypotheses apply to each incident edge, also when several
edges have the same image. The coefficient argument in the
introduction proves the domain assertion over every field.

For hyperbolicity, the edge images are quasiconvex on each side.
At a lift of a \(G^\pi\)-vertex, two distinct incident edges are
indexed by distinct pairs \((i,gF_i)\). Their stabilizers intersect
trivially by joint malnormality. Every three-edge segment has such
an interior vertex, so
\begin{equation}\label{eq:perm-acylindrical}
 \operatorname{Fix}_{\Gamma_p}(J)=1
 \quad\text{for every segment \(J\) of three edges}.
\end{equation}
This length matters. At the central vertex, two edges over the
same \(A_i\) and coset can have a common nontrivial stabilizer.

The vertex groups are hyperbolic, and all edge images are proper.
On the \(A\)-side \(N\ge2\) gives properness. On the other side,
\(F_i=G_p\) would contradict joint malnormality with another
nontrivial \(F_j\). The vertex and edge groups contain nonabelian
free groups. Apply Dahmani's combination theorem to this finite
acylindrical graph with empty vertex peripheral families. A proper
quasiconvex subgroup of a hyperbolic group is fully quasiconvex in
this setting \cite[Proposition 1.11]{r8}. The resulting group is
relatively hyperbolic with empty peripheral family, hence
word-hyperbolic \cite[Theorem 0.1(1)]{r8}, \cite[Theorem 0.1]{r5}.
Related combination theorems are given in \cite{r2,r3,r15}.
This completes the deduction.
\end{proof}

\section{Constructing the upper bound on rank}

What group relations can force a sum of ranks to be small? An outer product has rank at most one over every division ring. We impose triangle relations so that certain \emph{sums} of line matrices are outer products. Interpolation then recovers the individual line matrices from these sums. We describe the required configuration of lines first and prove its existence at the end of the section.

\subsection{Lines, points and blocks}

Fix a prime \(p\). We use the marked-line construction of \cite[Lemma~2.2]{r20} with larger girth and dimension. Set
\begin{equation}\label{eq:params}
 g=100,\quad d=g+3,\quad c=(20d)^{-d},\quad r=100c^{-1},\quad
 Q=p^a>Q_0=320000\cdot2060^{309}.
\end{equation}
Here \(g\) controls short cycles, \(d\) is the affine dimension, and \(r\) is the matrix size. None of these constants is optimized: each only has to be large compared with the ones before it, and Table~\ref{tab:parameters} records where each is used.

\begin{table}[t]
\centering\tabsmall
\renewcommand{\arraystretch}{1.12}
\begin{tabularx}{\linewidth}{@{}llX@{}}
\toprule
symbol & value & what it is for\\
\midrule
\(g\) & \(100\) & girth of the incidence graph; long circuits in the tip links make \(\widetilde X\) a \(\mathrm{CAT}(-1)\) space (Section~5)\\
\(d\) & \(g+3\) & dimension of the affine space \(P=\F_Q^d\) of points\\
\(c\) & \((20d)^{-d}\) & density of the lines, \(N\ge cm/2\) (Lemma~\ref{lem:lines})\\
\(r\) & \(100c^{-1}\) & size of the line matrices; every regular point carries at least \(4r^2\) labels\\
\(Q=p^a\) & \(>Q_0\) & size of the field; \(Q_0\) is the largest of \(20d\), \(32r^2/c\), \(256r/c\) and \(16rg^2\), as the proof of Lemma~\ref{lem:lines} requires\\
\(S\) & \(\lfloor Q/(10d)\rfloor\) elements & grid for the marked points; it keeps the interpolation degree \(2d(|S|-1)\) below \(Q-1\) in \eqref{eq:punctured}\\
\(\mathcal E\) & \(|\mathcal E|\le m/(2r)\) & exceptional points, where only the trivial rank bound \(r\) is used\\
\bottomrule
\end{tabularx}
\caption{The parameters of Section~4.}
\label{tab:parameters}
\end{table}

 Write \(P=\mathbb F_Q^d\) and \(m=|P|=Q^d\). Choose \(S\subset\mathbb F_Q\) of cardinality \(\lfloor Q/(10d)\rfloor\). On each line \(\ell_i\) mark a point \(a_i\in S^d\) and put \(\ell_i^\circ=\ell_i\setminus\{a_i\}\). The marks are distinct. The point--line incidence graph joins \(z\in P\) to the label \(i\) when \(z\in\ell_i^\circ\). Write \(I(z)\) for the set of labels at \(z\).

\begin{lemma}[Marked lines]\label{lem:lines}
There are \(N\ge cm/2\) marked lines and a set \(\mathcal E\subset P\) of size at most \(m/(2r)\) such that the incidence graph has girth greater than \(g\) and \(|I(z)|\ge4r^2\) for every \(z\notin\mathcal E\). In particular,
\begin{equation}
 2m+r|\mathcal E|\le5m/2<25m/2\le Nr/4.
\end{equation}
\end{lemma}

The points outside \(\mathcal E\) will each cost at most two units of rank. At each such point partition \(I(z)\) into two rectangles. To see that this is possible for every \(t=|I(z)|\ge4r^2\), set
\[
 b_t=r+(t\bmod r),\qquad a_t=(t-(r+1)b_t)/r.
\]
Then \(r\le b_t\le2r-1\) and \(a_t\) is an integer at least \(r\). The two rectangles of sizes \(r\times a_t\) and \((r+1)\times b_t\) have a total of \(t\) positions. Order the labels, use the initial segment for the first rectangle, and in each rectangle list the labels in lexicographic order. We write \(b\) for one of these blocks, \(z_b\) for its point, and \(I_b\) for its labels. Its coordinates are a bijection
\begin{equation}
 i\longmapsto(\alpha_b(i),\beta_b(i)):
 I_b\longrightarrow\mathbb Z/s_b\mathbb Z\times\mathbb Z/t_b\mathbb Z,\qquad s_b,t_b\ge r .
\end{equation}
There are no blocks at an exceptional point. The bipartite graph of blocks and their labels has girth greater than \(g\) as well. The map \(b\mapsto z_b\) to the incidence graph is locally injective, because a label belongs to at most one block at any given point.

\subsection{Triangles that give an outer product}

Write \([r]=\{0,\ldots,r-1\}\). The vertices of the complex \(X\) are the apices \(O_b\) and two sets of tips \(V_u,W_v\) indexed by \([r]\). There is an edge \(x_{iuv}:V_u\to W_v\) for each label \(i\) and each pair \(u,v\). From \(O_b\) to \(V_u\) draw edges \(L_{b,u,j}\) indexed by \(j\in\mathbb Z/s_b\mathbb Z\), and from \(O_b\) to \(W_v\) edges \(R_{b,v,k}\) indexed by \(k\in\mathbb Z/t_b\mathbb Z\). For \(i\in I_b\) and \(u,v\in[r]\) attach the triangle with oriented boundary
\begin{equation}\label{eq:triangle}
 L_{b,u,\alpha_b(i)+v}\,x_{iuv}\,R_{b,v,\beta_b(i)+u}^{-1}.
\end{equation}
Unused edges remain in \(X\). Each triangle is embedded, and the complex is finite and connected. For instance, a regular block joins all \(V_u,W_v\), and each apex is attached to these tips. Put \(G_p=\pi_1X\). Choose paths from a common basepoint to every vertex. For an oriented edge \(e:T\to T'\), the based element is \(g_e=[p_Tep_{T'}^{-1}]\). Define
\begin{equation}
 Y_i=(g_{x_{iuv}})_{u,v\in[r]}.
\end{equation}
These group-element matrices can be read over any coefficient field.

\begin{figure}[t]
\centering
\begin{tikzpicture}[font=\small, scale=1.05, >=Stealth]
\coordinate (O) at (0,2.25);
\coordinate (V) at (-2.1,0);
\coordinate (W) at (2.1,0);
\fill[black!7] (O)--(V)--(W)--cycle;
\draw[thick,->] (O)--node[pos=.55,left=4pt]{\(L_{b,u,\alpha_b(i)+v}\)} (V);
\draw[thick,->] (V)--node[below=2pt]{\(x_{iuv}\)} (W);
\draw[thick,->] (O)--node[pos=.55,right=4pt]{\(R_{b,v,\beta_b(i)+u}\)} (W);
\draw[dashed,gray] (O)--(-3.7,1.95) node[left,black]{\(V_{u'}\)};
\draw[dashed,gray] (O)--(3.7,1.95) node[right,black]{\(W_{v'}\)};
\fill (O) circle (1.6pt) node[above=2pt]{\(O_b\)};
\fill (V) circle (1.6pt) node[below left]{\(V_u\)};
\fill (W) circle (1.6pt) node[below right]{\(W_v\)};
\node at (0,0.75) {label \(i\in I_b\)};
\end{tikzpicture}
\caption{One triangle of the block \(b\) at the point \(z_b\), with the boundary \eqref{eq:triangle}. In the group it gives \(g_{x_{iuv}}=g_L^{-1}g_R\), where the \(L\)-edge depends only on \(u\) and the index \(\alpha_b(i)+v\), and the \(R\)-edge only on \(v\) and \(\beta_b(i)+u\). The apex \(O_b\) is joined in the same way to every tip \(V_{u'}\) and \(W_{v'}\) (dashed). As \(i\) runs through \(I_b\) the indices run through full residue systems, which gives the outer product \eqref{eq:outer}.}
\label{fig:triangle}
\end{figure}
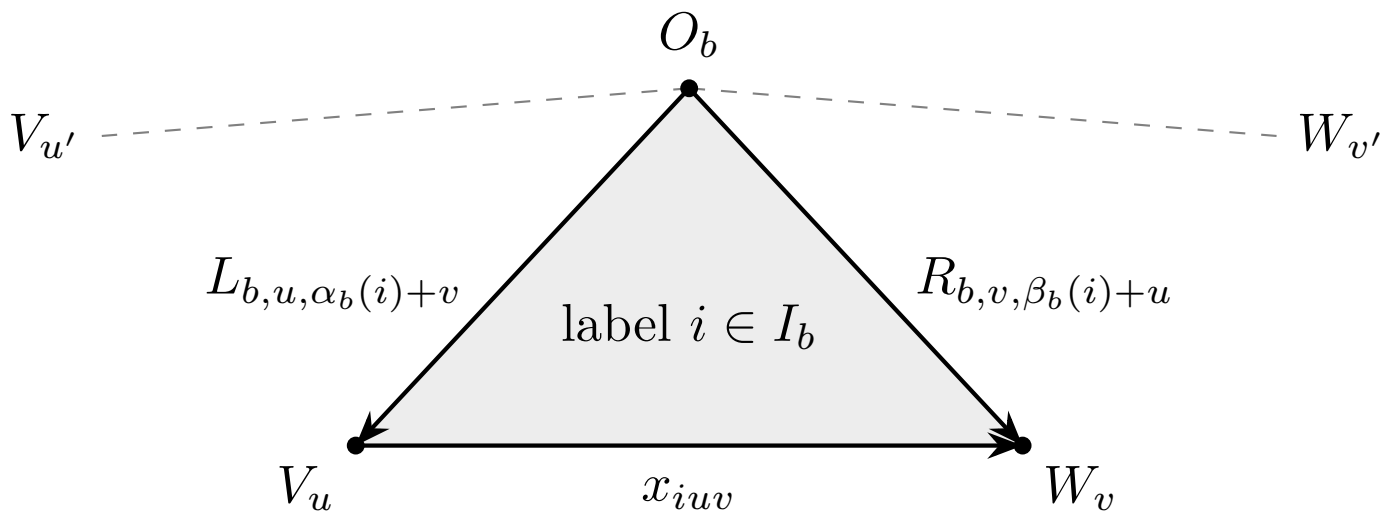

The reason for the shifts \(+v,+u\) is already visible in a single block. Suppressing a homomorphism to a division ring in the notation, (4.4) gives \(g_{x_{iuv}}=g_{L_{b,u,\alpha_b(i)+v}}^{-1}g_{R_{b,v,\beta_b(i)+u}}\). Summing over \(i\in I_b\) runs through a full rectangle, even after the two shifts. Thus
\begin{equation}\label{eq:outer}
 \left(\sum_{i\in I_b}Y_i\right)_{uv}
 =\left(\sum_j g_{L_{b,u,j}}^{-1}\right)\left(\sum_k g_{R_{b,v,k}}\right).
\end{equation}
The first factor depends only on \(u\), the second only on \(v\). This is the promised outer product; Figure~\ref{fig:triangle} shows the triangle that produces one of its terms. The shifts will also keep paths in the vertex links from folding in Section~5.

\subsection{Recovering the individual lines by interpolation}

Let \(I\) be the surviving label set. For \(i\in I\), take the Lagrange polynomial on the grid \(S^d\) that singles out \(a_i\):
\[
 f_i(X)=\prod_{\nu=1}^d\prod_{\sigma\in S,\ \sigma\ne a_{i,\nu}}
 \frac{X_\nu-\sigma}{a_{i,\nu}-\sigma}.
\]
In the column \(f=(f_i)_{i\in I}\), evaluation at \(a_i\) gives the \(i\)-th coordinate vector. Each entry of \(ff^T\) has total degree at most \(2d(|S|-1)<Q-1\). On a full line every such entry sums to zero. For constants this is \(Q=0\) in \(\mathbb F_Q\). For a monomial of degree \(1\le j<Q-1\), scale the line parameter by a nonzero \(\omega\) with \(\omega^j\ne1\), which exists by the root bound. Removing the marked point from the sum therefore leaves
\begin{equation}\label{eq:punctured}
 \sum_{z\in\ell_i^\circ} f(z)f(z)^T=-E_{ii}.
\end{equation}

\begin{lemma}[Interpolation and rank]\label{lem:interp}
Let \(D\) be any division ring of characteristic \(p\) and \(Y_i\) any \(r\times r\) matrices over \(D\). For \(Z(z)=\sum_{i\in I(z)}Y_i\),
\begin{equation}
 \sum_{i\in I}\operatorname{rank}_D Y_i\le\sum_{z\in P}\operatorname{rank}_D Z(z).
\end{equation}
\end{lemma}

\begin{proof}
We explain how to use the scalars in (4.7) when \(D\) contains only the specified prime field. Represent \(\lambda(c)\), multiplication by \(c\in\mathbb F_Q\), as an \(a\times a\) matrix over \(\mathbb F_p\). In \(M_a(D)\), each \(\lambda(c)\) commutes with each \(yI_a\), \(y\in D\), because the prime field is central. Thus \(\mathbb F_Q\) acts centrally in the subring generated by these two images. Read \(f(z)\) and each \(Y_i\) in that subring. Interchanging the two finite sums and using (4.7) gives
\begin{equation}
 -\bigoplus_iY_i=\sum_{z\in P}(f(z)\otimes I_r)\,Z(z)\,(f(z)^T\otimes I_r).
\end{equation}
Indeed, the coefficient of \(Y_i\) on the right is the block matrix \(-E_{ii}\). Now expand all \(a\times a\) blocks as matrices over \(D\). Every original matrix \(Y\) becomes \(Y\otimes I_a\), which is a direct sum of \(a\) copies of \(Y\) after reordering rows and columns. The left side has rank \(a\sum_i\operatorname{rank}_D Y_i\), while the summand at \(z\) has rank at most \(a\operatorname{rank}_D Z(z)\). Subadditivity and cancellation of \(a\) give the result.
\end{proof}

\begin{proposition}[The rank bound]
Every homomorphism \(\phi:\mathbb F_p[G_p]\to D\) to a division ring satisfies
\begin{equation}
 \sum_i\operatorname{rank}_D\phi(Y_i)\le2(m-|\mathcal E|)+r|\mathcal E|<Nr/4.
\end{equation}
\end{proposition}

\begin{proof}
At each regular point, \(\phi(Z(z))\) is the sum of the two outer products in (4.6) and has rank at most two. At an exceptional point the matrix size gives the bound \(r\). Apply Lemma~4.2 and (4.2). Unitality ensures that \(D\) has characteristic \(p\).
\end{proof}

\subsection{Existence of the incidence configuration}

\begin{proof}[Proof of Lemma~4.1]
The threshold \(Q_0\) in (4.1) is the largest of \(20d\), \(32r^2/c\), \(256r/c\) and \(16rg^2\). At each of the \(n_0=|S|^d\) marks choose one line direction independently and uniformly out of \(J=(Q^d-1)/(Q-1)\) possibilities. Initially two marks may choose the same geometric line. We have \(n_0\ge cm\) and \(Q^{d-1}\le J\le2Q^{d-1}\).

First count points of low degree. For fixed \(z\), each mark other than \(z\) contributes a punctured incidence with probability \(1/J\). Their sum is binomial with variance at most its mean \(\mu_z\ge(n_0-1)/(2Q^{d-1})\ge cQ/4>8r^2\). By Chebyshev's inequality the probability of degree below \(4r^2\) is at most \(4/\mu_z\le16/(cQ)\). Hence, by Markov's inequality, the probability of more than \(m/(4r)\) low-degree points is at most \(64r/(cQ)<1/4\).

Next count short cycles. A simple cycle of length \(2s\), for \(2\le s\le g/2\), uses at most \(m^s\) ordered choices of point vertices. Consecutive points determine a line on which at most \(Q\) marks are possible. For a valid choice of distinct marks, the required directions have probability \(J^{-s}\). Thus its length class has expected count at most \(m^s Q^s J^{-s}\le Q^{2s}\). The expected number of all short cycles is at most \(gQ^g\). With probability at least \(3/4\) there are at most \(4gQ^g\) such cycles. Choose an outcome meeting both bounds.

Delete every label on a short cycle. At most \(4g^2Q^g\) labels and at most \(4g^2Q^{g+1}\) points on their lines are involved. Since \(d=g+3\), these counts divided by \(m\) are at most \(4g^2/Q^3\) and \(4g^2/Q^2\), respectively. They are bounded by \(c/2\) and \(1/(4r)\) because \(Q>16rg^2\) and \(r=100/c\). Take \(\mathcal E\) to consist of the original low-degree points and the points on deleted lines. All other degrees are unchanged, the surviving graph has the required girth, and \(N\ge cm/2\), \(|\mathcal E|\le m/(2r)\). Finally \(Nr/4\ge cmr/8=25m/2\), giving (4.2).
\end{proof}

All choices can be made by finite search. Section~7 specifies the orders that select a particular presentation.

\section{The geometry of the base group}

The triangle relations have already supplied the upper bound on rank. We now need the entries of each normalized matrix to be free generators, and we need \(K[G_p]\) to be a domain. The first assertion is detected in the links of vertices. The second requires a length estimate for essential loops in the whole complex. The outputs are Lemma~\ref{lem:freebasis}, which gives Proposition~\ref{prop:base}(i), the quasiconvexity and joint malnormality used in Section~3, and Corollary~\ref{cor:diffuse}, the diffuseness of \(G_p\).

\subsection{A metric on the triangles}

Give every triangle of \(X\) the hyperbolic right isosceles metric with right angle at \(O_b\) and legs of length \(3\). Write \(a\) for its hypotenuse length, \(\theta\) for a tip angle and \(h\) for the altitude to the hypotenuse. Give unused \(x\)-edges length \(a\). The right-triangle identities give
\begin{equation}
 \cosh a=\cosh^2 3,\qquad
 \sin\theta=(\cosh^2 3+1)^{-1/2},\qquad
 \sinh^2h=\frac{\cosh^2 3-1}{\cosh^2 3+1}.
\end{equation}
At an apex the link is a complete bipartite graph with edges of length \(\pi/2\). At a tip the links can be more complicated, but their geometry is controlled by the incidence graph.

\begin{lemma}
The universal cover \(\widetilde X\) is proper \(\mathrm{CAT}(-1)\). The group \(G_p\) is torsion-free and word-hyperbolic.
\end{lemma}

\begin{proof}
An \(L_{b,u,j}\)-germ and an \(R_{b,v,k}\)-germ at \(O_b\) determine one triangle, since its label must have coordinates \((j-v,k-u)\). Every circuit in the apex link therefore has length at least \(2\pi\). In the link at \(V_u\) send \(L_{b,u,j}\) to \(z_b\) and \(x_{iuv}\) to \(i\). This graph map is locally injective. At an \(L\)-germ the condition \(\alpha_b(i)+v=j\) determines at most one \(v\in[r]\), as \(s_b\ge r\). At an \(x\)-germ there is at most one incident block over a given point. At \(W_v\) the same proof uses \(t_b\ge r\). Each reduced circuit in a tip link therefore has more than \(g\) edges. Since \(10<\cosh 3<11\), (5.1) gives \(\theta>\sin\theta>1/12\), and \(g\theta>2\pi\).

The link condition and Cartan--Hadamard give the asserted curvature bound \cite[I.7.50, II.5.2, II.5.6 and II.4.1(2)]{r7}. Local finiteness and finite shapes give properness. The deck action is free, proper and cocompact. A finite-order isometry would fix a point \cite[II.2.8(1)]{r7}. The orbit is quasi-isometric to \(G_p\) by \cite[I.8.19]{r7}, so \(G_p\) is word-hyperbolic.
\end{proof}

\subsection{The free subgroups in the matrices}

The edges \(x_{iuv}\) of one label form a copy \(\Theta_i\) of the complete bipartite graph \(K_{r,r}\). Its row zero and column zero span a tree. Each remaining edge specifies a basis loop, and those loops are precisely the entries that will survive normalization.

\begin{lemma}\label{lem:freebasis}
The map \(\pi_1\Theta_i\to\pi_1X\) is injective. For \(g_{iuv}=g_{x_{iuv}}\), its image \(F_i\), based by the path to \(W_0\), has free basis
\begin{equation}
 e_{i,uv}=g_{iu0}^{-1}g_{iuv}g_{i0v}^{-1}g_{i00}\qquad(1\le u,v<r).
\end{equation}
\end{lemma}

\begin{proof}
Distinct germs of label \(i\) at a tip are far apart in the link. Any reduced path between them projects by local injectivity to a nonempty incidence walk from \(i\) back to \(i\), with no internal immediate reversal. A shortest positive segment with coinciding endpoints is a simple cycle, so the path has more than \(g\) edges. Thus their link distance is greater than \(g\theta>2\pi\), or is infinite. Reduced paths in \(\Theta_i\) are consequently locally geodesic in \(X\). In the universal cover such paths are geodesics, so a nontrivial reduced loop cannot contract. The non-tree edges with the chosen paths through row and column zero give the displayed basis.
\end{proof}

Multiplying \(Y_i\) by the diagonal matrices with entries \(g_{iu0}^{-1}\) on the left and \(g_{i0v}^{-1}g_{i00}\) on the right gives precisely the first matrix in (2.1). Subtract row zero from each other row, then clear the other entries of row zero by column operations. This proves Proposition~2.1(i).

\begin{lemma}
Each \(F_i\) is quasiconvex, and
\begin{equation}
 F_i\cap gF_jg^{-1}\ne 1\qquad\Longrightarrow\qquad i=j,\quad g\in F_i .
\end{equation}
\end{lemma}

\begin{proof}
A lifted line graph is a tree \(T\). By the link estimate just proved, every tree segment is a geodesic of \(\widetilde X\). Uniqueness of geodesics makes \(T\) convex. Its stabilizer, a conjugate of \(F_i\), acts cocompactly on \(T\), so that its orbit is quasiconvex there and in \(\widetilde X\). The orbit quasi-isometry and stability of quasigeodesics give quasiconvexity in a Cayley graph of \(G_p\) \cite[I.8.19, III.H.1.7]{r7}.

A nonidentity deck transformation stabilizing two such trees would translate along an axis in each. These two axes are asymptotic and bound a flat strip \cite[II.2.13]{r7}. In curvature at most \(-1\) the strip has width zero, so the trees share an entire line. Different lifts of the same graph are disjoint, and lifts of different labels share no edge. Thus the trees and their labels coincide. The coset indexing of lifts gives (5.3).
\end{proof}

\subsection{From lengths to diffuseness}

To prove diffuseness it is enough to find two extremal points in each finite orbit subset. We take the endpoints of a diameter. The following comparison ensures that at least one of two opposite translates of either endpoint moves farther from the other. The same criterion will apply to a quotient with shorter essential loops in Section~6.

\begin{lemma}[Midpoint criterion]
Let a group act freely by isometries on a \(\mathrm{CAT}(-1)\) space. If every nonidentity element has translation length at least \(\operatorname{arcosh}3\), then the group is diffuse.
\end{lemma}

\begin{proof}
Write \(d=d(x,\gamma x)\), \(e=d(x,\gamma^2 x)\) for \(\gamma\ne1\). If \(z\) is the midpoint of \([x,\gamma x]\), then \(\gamma z\) is the midpoint of the next side of the isosceles triangle. Write \(L=\inf_y d(y,\gamma y)\). The hyperbolic comparison triangle has equal sides \(d\). Applying the cosine law at their common vertex, first at distance \(d\) and then at distance \(d/2\), gives
\[
 \cosh L\le\cosh d(z,\gamma z)\le1+\frac{\cosh e-1}{2(\cosh d+1)}.
\]
Indeed, the ratio of the two increments above \(1\) is \(\sinh^2(d/2)/\sinh^2d=1/[2(\cosh d+1)]\). As \(\cosh L\ge3\), it follows that \(\cosh e\ge4\cosh d+5\). Since \(\cosh(d+\log 4)<4\cosh d\), we obtain \(e>d+\log 4\).

The strong four-point inequality in a \(\mathrm{CAT}(-1)\) space \cite[Definition~4.1 and Theorem~5.1]{r19} is
\begin{equation}
 e^{(d(a,b)+d(c,d'))/2}\le e^{(d(a,c)+d(b,d'))/2}+e^{(d(a,d')+d(b,c))/2}.
\end{equation}
Apply it to \(o,x,\gamma x,\gamma^{-1}x\). If both translates were at most \(d(o,x)\) from \(o\), it would give \(e^{(d(o,x)+e)/2}\le2e^{(d(o,x)+d)/2}\), contradicting the strict inequality above. Now choose any two diameter endpoints of a finite subset of a free orbit. With either endpoint as \(o\), the other is extremal. They give the two required points.
\end{proof}

For hyperbolic manifolds Bowditch proved diffuseness when the length bound is strictly greater than \(\operatorname{arcosh}3\) \cite[Theorem~5.3 and Lemma~5.4]{r6}. In particular, a free, proper, cocompact action on a proper \(\mathrm{CAT}(-1)\) space with all nontrivial translation lengths at least \(3\) satisfies the criterion. Our complex has this latter bound. The individual cone fillings will have the bound \(2\), still above \(\operatorname{arcosh}3\).

\subsection{How a short loop would travel through the complex}

Recall the altitude \(h\) in (5.1). For the length estimates put
\begin{equation}
 \eta=\pi/4-\arctan\!\left(\frac{\tanh(3/2)}{\sinh h}\right),\qquad
 \varepsilon=2\operatorname{arsinh}(\sinh h\sin\eta).
\end{equation}
The only numerical margins used in the argument are
\begin{equation}
 h>3/4,\qquad \eta>0,\qquad \varepsilon>79/1000,\qquad
 g\theta>2\pi,\qquad 6/\varepsilon<g .
\end{equation}
For completeness, \(20<e^3<21\) and \(e^{3/4}<9/4\) by the exponential series. Hence \(A=\sinh h>99/100>\sinh(3/4)\) and \(B=\tanh(3/2)<10/11\), with \(A,B<1\). Thus
\[
 A\sin\eta=\frac{A(A-B)}{\sqrt{2(A^2+B^2)}}
 >\frac{99}{100}\frac{99/100-10/11}{2}=\frac{8811}{220000}>\frac{1}{25}.
\]
Integration of \(1/\sqrt{1+t^2}\) from \(0\) to \(1/25\) yields \(\varepsilon>2/\sqrt{626}>79/1000\). The bound for \(g\theta\) was proved in Lemma~5.1.

\begin{figure}[tbp]
\centering
\begin{tikzpicture}[scale=3.05,>=Stealth,font=\small]
 % Poincare coordinates for a right isosceles triangle of side length 3.
 \pgfmathsetmacro{\rad}{(exp(3)-1)/(exp(3)+1)}
 \pgfmathsetmacro{\tip}{\rad/sqrt(2)}
 \pgfmathsetmacro{\cen}{(1+\rad*\rad)/(2*\tip)}
 \pgfmathsetmacro{\cir}{sqrt(\cen*\cen-1)}
 \pgfmathsetmacro{\ang}{atan2(\cen-\tip,\tip)}
 \pgfmathsetmacro{\foot}{\cen-\cir}
 \begin{scope}[xshift=-1.12cm]
  \coordinate (O) at (0,0);
  \coordinate (V) at (-\tip,-\tip);
  \coordinate (W) at (\tip,-\tip);
  \draw[semithick] (V)--(O)--(W);
  \draw[semithick] (W) arc[start angle=\ang,end angle=180-\ang,radius=\cir];
  \draw[densely dashed] (O)--(0,-\foot) node[midway,right] {$h$};
  \draw (-.045,-.045)--(0,-.09)--(.045,-.045);
  \fill (O) circle (.012) node[above] {$O_b$};
  \fill (V) circle (.012) node[below left] {$V_u$};
  \fill (W) circle (.012) node[below right] {$W_v$};
  \node[left] at (-.33,-.29) {$3$};
  \node[right] at (.33,-.29) {$3$};
  \node[below] at (0,-\foot-.025) {$x_{iuv}$};
 \end{scope}
 \begin{scope}[xshift=1.12cm,yshift=-.29cm]
  \draw[gray,densely dashed] (0,-.78)--(0,.78);
  \draw[semithick] (\tip,-\tip)
    arc[start angle=\ang,end angle=180-\ang,radius=\cir];
  \draw[semithick] (\tip,\tip)
    arc[start angle=-\ang,end angle=-180+\ang,radius=\cir];
  \draw[thick] (0,-\foot)--(0,\foot);
  \fill (0,0) circle (.012) node[right] {$O_b$};
  \node[right] at (0,\foot/2) {$h$};
  \node[right] at (0,-\foot/2) {$h$};
  \draw (0,\foot-.055)--(.055,\foot-.055)--(.055,\foot);
  \draw (0,-\foot+.055)--(.055,-\foot+.055)--(.055,-\foot);
 \end{scope}
\end{tikzpicture}
\caption{A right isosceles cell in the Poincar\'e disk (left). Unfolding the link sectors for two edges of one line places their hypotenuses on opposite sides of the apex (right). Their common perpendicular has length \(2h\).}
\label{fig:cell}
\end{figure}
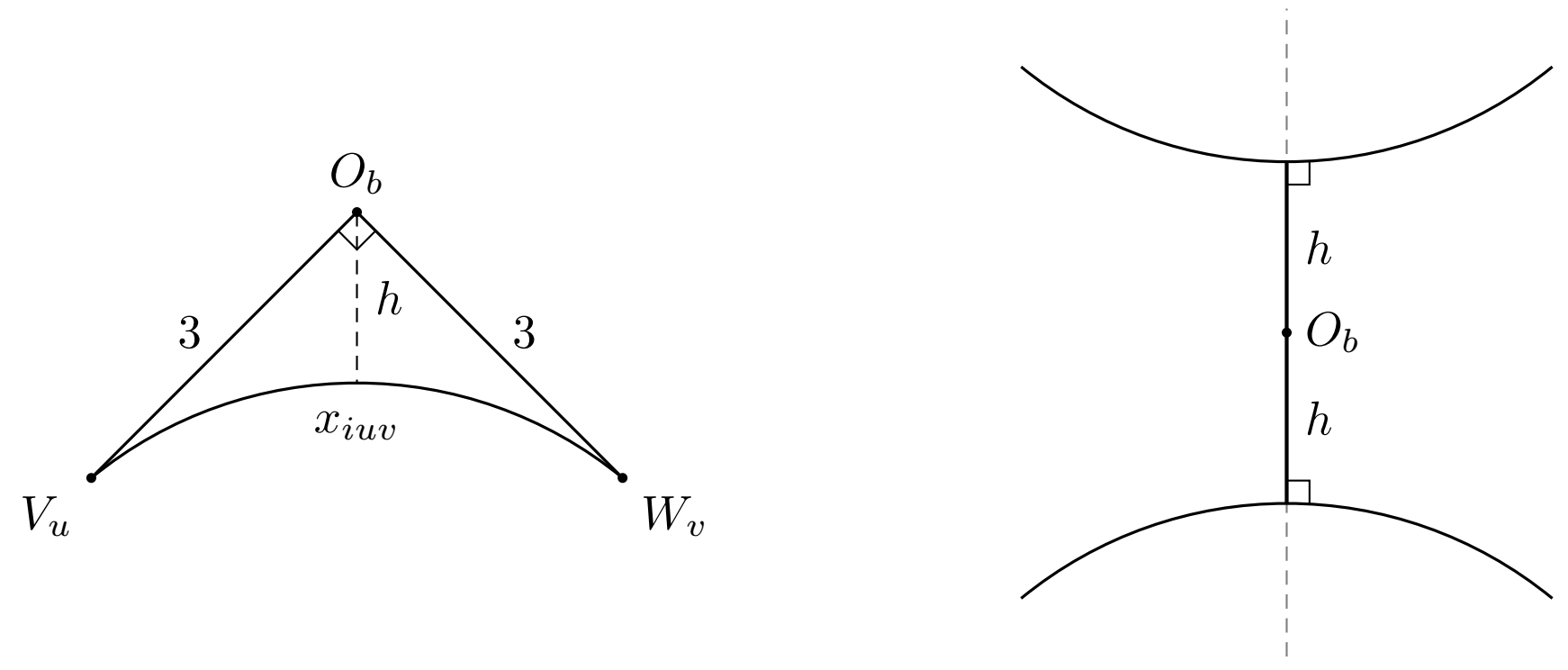

\begin{proposition}
Every essential loop in \(X\) has length at least \(3\).
\end{proposition}

\begin{proof}
Suppose there is a shorter essential loop. We may take it to be a finite string of cell segments. Indeed its lift joins a point to a nontrivial deck translate at distance less than \(3\), and the polyhedral length metric supplies such a string with length still less than \(3\).

\emph{Where the loop can go.} The open star of a tip is contractible and contains the open radius-\(3\) ball about it. In an incident triangle, the distance from the tip to the opposite radial side is at least \(3\), as the perpendicular meets that side at the apex. Unused incident edges have length \(a\ge3\). These radial distance functions agree on shared sides and are \(1\)-Lipschitz up to the first exit from the star. If the loop came within \(3/2\) of a tip, its shorter arcs would put it entirely within that star. Thus every crossing of an \(x\)-edge is at least \(3/2\) along the edge from each endpoint. In signed distance \(s\) from the midpoint this gives
\begin{equation}
 |s|\le a/2-3/2\le3/2.
\end{equation}

With the tips removed, the triangles and radial edges at one apex form a contractible block, including its open hypotenuses. Distinct blocks meet only along open \(x\)-edges. An unused open edge alone cannot carry the essential loop. Record the block passages in cyclic order. Merge successive passages in the same block. A passage that returns through the same \(x\)-edge may be replaced by the segment of that edge and assigned to the adjacent block. This replacement is a homotopy in the contractible block and does not increase length, since the lifted edge segment is geodesic. It stays in the window (5.7). Each reduction removes a passage. A record supported on one block contracts. We are left with successive distinct blocks and, inside each block, two distinct crossing edges.

\emph{The cost of one passage.} Join the two crossing points to the apex in a lift of their block. For a point with midpoint coordinate \(s\), the angle from the midpoint direction has tangent \(\tanh|s|/\sinh h\), by the right-triangle formulas. Its direction is therefore at least \(\eta\) from both ends of its edge in the apex link. Two distinct hypotenuses correspond to distinct link edges, so the angle between the directions is at least \(2\eta\). Their radial distances \(R_1,R_2\) are at least \(h\). Hinge comparison gives
\[
 \cosh D\ge\cosh(R_1-R_2)+2\sinh R_1\sinh R_2\sin^2\eta
 \ge1+2\sinh^2h\sin^2\eta.
\]
Thus each passage costs at least \(\varepsilon\).

A passage between two edges of the same line costs at least \(2h\). Here the two link edges have disjoint endpoints. Indeed, if the row is the same, different columns give different \(L\)-indices modulo \(s_b\) and different \(R\)-tips. The column case is symmetric, and if both indices differ both tips differ. If the radial directions have link distance at least \(\pi\), their geodesic goes through the apex. Otherwise a shortest link path traverses one whole intervening edge, since the edges have disjoint endpoints and two whole edges already have length \(\pi\). Unfold the three sectors in the plane of curvature \(-1\). The two midpoint directions are opposite, their angular separation being \(\pi/4+\pi/2+\pi/4=\pi\). The supporting lines of the hypotenuses have a common perpendicular of length \(2h\) (Figure~\ref{fig:cell}). Hence any two points on them are at least \(2h\) apart. Hinge comparison gives the same lower bound in the universal cover.

\emph{Why there must be two expensive passages.} If \(k\) passages remain, then \(k\varepsilon<3\). Forget the row and column of each crossing edge, keeping its line label. The resulting closed walk in the block--line graph has length \(2k<6/\varepsilon<g\). Its image is a finite tree with at least one edge. At each leaf the walk must reverse. A line cannot be such a leaf, since that would repeat consecutive blocks. Each of at least two block leaves therefore corresponds to leaving on a different edge of the same line. These passages alone cost at least \(4h>3\), a contradiction.
\end{proof}

\begin{corollary}\label{cor:diffuse}
The base group \(G_p\) is diffuse. Its group algebra is a domain over every field.
\end{corollary}

\begin{proof}
A path from a point in \(\widetilde X\) to a nontrivial deck translate projects to an essential loop. Thus every nontrivial translation length is at least \(3\). Apply Lemma~5.4 and the unique-product argument in the introduction.
\end{proof}

\section{Why the line subgroups are concentrated}

Quasiconvexity and malnormality gave hyperbolicity of the permutation group.
For its group algebra to be a domain we also needed diffuseness of the action
on $G_p/F_i$. We prove that property one line at a time. Fix $i$, let $L_i$
be the normal closure of $F_i$ in $G_p$, and write $Q_i=G_p/L_i$. Two facts
will suffice. The quotient $Q_i$ is diffuse, and $F_i$ is a free factor of the
free group $L_i$. The first fact follows by filling the line graph with a cone.
For the second we lift that cone and remove its apex from every sheet. Only the hyperbolic group of Theorem~\ref{thm:A} uses this section; Theorems~\ref{thm:B}--\ref{thm:D} need only Sections~4 and~5.

\subsection{The quotient by one line}

Attach the cone on $\Theta_i$ to $X$, using radial edges of length $3$ from a
new apex and the same right isosceles triangles along $x_{iuv}$. Call the
resulting finite complex $Z_i$. Van Kampen gives $\pi_1 Z_i=Q_i$. These are
auxiliary spaces, one for each $i$. No cone is attached in the definition of
$\Gamma_p$.

\begin{lemma}
The universal cover of $Z_i$ is proper $\mathrm{CAT}(-1)$, its quotient $Z_i$
has no essential loop of length less than $2$, and $Q_i$ is diffuse.
\end{lemma}

\begin{proof}
The link at the new apex is $K_{r,r}$ with edge lengths $\pi/2$. At an old
tip, a new circuit uses the new radial vertex and a path in the old link
between distinct germs of label $i$. That path has length greater than
$g\theta$ by the proof of Lemma~5.2. All vertex links still have circuits of
length at least $2\pi$. Edge-interior links are unions of semicircles between
the tangent directions. The link condition and Cartan--Hadamard therefore
give the same proper $\mathrm{CAT}(-1)$ universal-cover conclusion as for
$X$ \cite[I.7.50, II.5.2, II.5.6 and II.4.1(2)]{r7}.

Suppose an essential loop has length less than $2$. A point within distance
$2$ of a tip would put the entire loop in its contractible open star, since
shorter arcs would put it within distance $3$ of that tip. Thus every
crossing point lies at least $2$ along its edge from both endpoints and
within $a/2-2\leq1$ of the midpoint. Apply the block reduction of
Proposition~5.5, now counting the new cone as a block $c$. The reduction
uses only contractibility of each tip-deleted block, geodesicity of the lifted
edges and convexity of the permitted interval on each edge. All three
properties remain valid.

Every old passage still costs at least $\varepsilon>79/1000$, or at least
$2h>3/2$ if it joins edges of the same line. A passage in $c$ costs more
than $1/2$. To prove this, first suppose its two edges share a tip. They meet
there at angle $2\theta$, since the new two-edge path in the tip link is
shorter than any old path between the two germs. If the two distances from
the tip are $s,t\geq2$, the cosine law and hinge comparison give directly
\begin{equation}
 \cosh D\ \geq\ \cosh(s-t)+2\sinh s\sinh t\sin^2\theta
       \ \geq\ 1+2\sinh^2(2)\sin^2\theta .
 \label{eq:cone-hinge}
\end{equation}
Thus $D\geq2\operatorname{arsinh}(\sinh 2\sin\theta)>1/2$, since
$\sinh 2>7/2$, $\sin\theta>1/12$, and $\sinh(1/4)<4/15<7/24$.
If the two edges share no tip, their radial directions at the cone apex have
link distance at least $\pi/2$, and both radial lengths are at least $h$.
This time comparison gives $\cosh D\geq\cosh^2 h>\cosh(1/2)$, using
$h>3/4$.

For $k$ passages the projected walk has $2k<4/(79/1000)<100$. Its graph is
the old block--line graph with a single leaf $c$ attached at $i$. Its image
must again be a tree with at least two leaves, all of them block vertices.
Two old leaves cost at least $4h>3$. If one is $c$, the other is old, and
together they cost more than $1/2+2h>2$. This excludes the loop. Every
nontrivial deck displacement is at least $2>\operatorname{arcosh}3$, so
Lemma~5.4 applies.
\end{proof}

\subsection{Removing the lifts of the apex}

Removing one point from a small cone leaves a copy of its link up to
homotopy. To identify the normal closure of $F_i$ we need to know how the
fundamental groups of all the links fit together. The following lemma
supplies compatible connecting paths. This choice matters when the factors
are viewed as subgroups.

\begin{lemma}[Puncturing a polyhedral complex]
Let $M$ be a proper, locally finite $\mathrm{CAT}(0)$ $M_\kappa$-polyhedral
complex, with $\kappa\leq0$ and finitely many shapes. Let $P$ be a locally
finite set of vertices with nonempty, connected links. Then $M\setminus P$
is path-connected. One can choose paths $\gamma_v$ from a common basepoint
to small link spheres about $v\in P$ so that inclusion gives an isomorphism
\begin{equation}
 \mathop{\ast}_{v\in P}\pi_1(\operatorname{Lk}_M v)
    \ \longrightarrow\ \pi_1(M\setminus P),\qquad
    [\alpha]\longmapsto[\gamma_v\alpha\gamma_v^{-1}].
 \label{eq:punctured-free-factors}
\end{equation}
\end{lemma}

\begin{proof}
If $P$ is empty, use contractibility of $M$. Otherwise fix $p\in P$. Grow a
closed ball about $p$, writing $Y_R=\overline B(p,R)\setminus P$. In a
sufficiently small punctured ball we see only $\operatorname{Lk}p$. Between
successive radii containing punctures, radial projection along the unique
geodesics toward $p$ retracts the larger $Y_R$ onto the smaller. This
homotopy moves points only in the band of radii being removed. It therefore
suffices to understand what is added at one puncture radius. Properness and
local finiteness give only finitely many punctures in every bounded ball.

Let $v_1,\ldots,v_k$ be the punctures at distance $d>0$ from $p$. Choose
$\delta>0$ so that $d-\delta>0$, the interval $[d-3\delta,d+3\delta]$
contains no other puncture radius, and the $2\delta$-balls about the $v_j$
are disjoint cone neighborhoods containing no other punctures
\cite[I.7.39]{r7}. Put $W=Y_{d+\delta}$. The obstruction to radial
projection at $v_j$ is its shadow
\[
 S_j=\{x\in\overline B(p,d+\delta):v_j\in[p,x]\}.
\]
This set is closed by continuous dependence of geodesics on their endpoints.
It lies in $\overline B(v_j,\delta)$, by additivity of length along $[p,x]$.
Removing the shadows gives an open set $O=W\setminus\bigcup_j S_j$ which
retracts radially onto $Y_{d-\delta}$. Indeed, every point on the radial
segment from a point of $O$ toward the smaller ball is still outside the
shadows.

The other members of an open cover of $W$ are $V_j=W\cap B(v_j,2\delta)$.
Each retracts to a copy of its link at a fixed radius $0<\eta<\delta$. In
cone coordinates about $v_j$, each direction contains the initial interval
$(0,\delta]$ in $V_j$, by the triangle inequality. Convexity of distance to
$p$ makes the full permitted interval initial, so interpolation to radius
$\eta$ stays in $V_j$.

It remains to check how $O$ meets $V_j$. Let $\xi$ be the direction from
$v_j$ toward $p$. For $x\neq v_j$ in that cone neighborhood, with direction
$u$, one has
\begin{equation}
 v_j\in[p,x]\quad\Longleftrightarrow\quad
 d_{\operatorname{Lk}v_j}(\xi,u)\geq\pi .
 \label{eq:link-shadow}
\end{equation}
The cone angle is the link distance capped at $\pi$. At angle $\pi$ the two
geodesic segments concatenate to a local geodesic, hence to a geodesic in
$M$ \cite[II.1.4]{r7}. Conversely, a geodesic through $v_j$ has angle $\pi$.
It follows that $O\cap V_j$ retracts to the open radius-$\pi$ ball about
$\xi$ in the link. This ball is nonempty and contractible. In a
$\mathrm{CAT}(1)$ link, every point at distance less than $\pi$ is joined to
the center by a unique, continuously varying geodesic contained in the
ball \cite[II.1.4]{r7}.

The $V_j$'s are disjoint and path-connected, and each has path-connected,
simply connected overlap with $O$. Van Kampen therefore adjoins one free
factor $\pi_1(\operatorname{Lk}v_j)$ for each $j$ and preserves
path-connectedness. Paths to the new factors can be chosen through $O$,
keeping the earlier choices. Starting with the connected link of $p$,
iterate over all puncture radii. Every loop and every homotopy of a loop has
compact image and lies in some $Y_R$. Thus the successive free-factor
inclusions give \eqref{eq:punctured-free-factors} in the union.
\end{proof}

\subsection{The normal closure and the coset action}

In $M=\widetilde Z_i$, let $P$ be the set of lifts of the new apex. It is
indexed by $Q_i$, and every vertex in $P$ has link $\Theta_i$. The inverse
image of $X$ is connected since $G_p\to Q_i$ is onto. It is the cover of
$X$ corresponding to $L_i$. Each lifted cone is embedded and its punctured
cone retracts onto its base. These retractions agree on the inverse image of
$X$, and local finiteness makes them a deformation retraction of
$M\setminus P$ onto that cover. Under this retraction the link inclusions
are the inclusions of the lifted line graphs. Returning the basepoints along
the paths in Lemma~6.2 gives, for appropriate representatives $g_q$,
\begin{equation}
 L_i=\mathop{\ast}_{q\in Q_i}g_q F_i g_q^{-1}.
 \label{eq:line-closure}
\end{equation}
In particular $L_i$ is free. Conjugating the decomposition to put one chosen
factor in its original position shows that $F_i$ itself is a free factor of
$L_i$, since $L_i$ is normal in $G_p$.

\begin{proof}[Completion of Proposition~2.1]
A free factor is concentrated in a free group, as shown in Section~3. Thus
$F_i$ is concentrated in $L_i$. The action on $G_p/L_i$ is the regular
action of $Q_i$ through the quotient map. Its diffuseness, proved in
Lemma~6.1, says exactly that $L_i$ is concentrated in $G_p$. Transitivity
of concentration \cite[Theorem~1.3(1)]{r6} now proves that $F_i$ is
concentrated in $G_p$. This was the last assertion needed for the two
constructions.
\end{proof}

\section{Presentations and the embedding question for the base}

\subsection{Fixing all choices}

Here are deterministic presentations of the examples. Choose the least $a$
with $Q=p^a>Q_0$ in (4.1), and the first monic irreducible polynomial of
degree $a$ in lexicographic order on its coefficient residues. Use its
power basis for $\mathbb F_Q$. Number $\sum_{j=0}^{a-1}c_jX^j$ by
$\sum_{j=0}^{a-1}c_jp^j$, with $0\leq c_j<p$, and use the first
$\lfloor Q/(10d)\rfloor$ elements for $S$. Normalize each direction by
making its first nonzero coordinate $1$. Order marks and normalized
directions lexicographically by these numbers, and choose the first
direction tuple satisfying the two bounds in the probabilistic proof before
deletion. Delete every label on a short cycle and take the indicated
exceptional set. Lemma~4.1 proves termination.

Use (4.3) for the blocks, order the vertices and edges by their displayed
indices (taking the types in the order of their introduction), and take the
first spanning tree. In the resulting presentation of $G_p$, its tree
generators are $1$, opposite orientations are inverses, and the relators are
precisely (4.4). Compute each line word from (5.2). For the shear group add
two letters $s_i,\gamma_i$ for each line, with relations
\begin{equation}
 [s_i,\gamma_i]=[s_i,e_{i,uv}]=1,\qquad
 \gamma_i e_{i,uv}\gamma_i^{-1}=s_i^{\delta_{uv}}e_{i,uv}
     \qquad(1\leq u,v<r).
 \label{eq:line-shear-presentation}
\end{equation}
For $\Gamma_p$ take the first line as $i_*$ and use instead (3.2)--(3.3).

\subsection{What an embedding of the base algebra would have to do}

The base is already a domain over every field. It also rules out embeddings
satisfying the coset-independence property below. The question left open is
whether the base admits any division-ring embedding at all.

Given $K[G]\subset D$ and $C\leq G$, let $D_C$ be the division subring
generated by $K[C]$. A Linnell division ring is a division ring generated
by $K[G]$ with the following property. For every $C\leq G$, elements in
distinct right cosets of $C$ are left linearly independent over $D_C$
\cite[Section~2.3]{r14}, \cite[Definition~2.2]{r11}.

\begin{proposition}
For every field $K$ of characteristic $p$, the domain $K[G_p]$ has no
Linnell division ring.
\end{proposition}

\begin{proof}
Assume there is such an embedding. On one line, suppress its label and let
$M_k=(e_{uv}-1)_{1\leq u,v\leq k}$. The first pivot $e_{11}-1$ is nonzero.
If $M_{k-1}$ is invertible, let $C$ be generated by all the line basis
elements except $e_{kk}$. In the bordering matrix with last column $b$,
last row $c$ and new corner $e_{kk}-1$, the next pivot is
\[
 e_{kk}-1-cM_{k-1}^{-1}b .
\]
Here $b,c$ denote only the entries outside the corner. The inverse of
$M_{k-1}$ is over $D_C$ by elimination, so the pivot is $e_{kk}$ plus an
element of $D_C$. The free basis puts $1$ and $e_{kk}$ in distinct right
cosets of $C$, so coset independence makes this pivot nonzero. Induction
gives full rank on every line, contrary to (4.10) after restriction to
$\mathbb F_p[G_p]$.
\end{proof}

This is the same obstruction to a Linnell embedding that occurs for the
earlier group in \cite[Theorem~1.4 and Corollary~6.5]{r9}. For an arbitrary
embedding the argument stops at the pivot, since coset independence was an
extra hypothesis. The shears force every pivot matrix to be invertible by
conjugation. The permutation graph only needs enough successful pivots on
average. This is the role of enlarging the base.

In particular $\Gamma_p$ and $\Lambda_p$ are counterexamples in characteristic $p$ to \cite[Conjecture~1]{r14} and \cite[Conjecture~1.1]{r11}. If $\mathbb F_p[G_p]$ does embed, it will give a domain embeddable in a
division ring but admitting no Linnell division ring. If it does not, $G_p$
itself will give a smaller hyperbolic example. Neither construction
determines embeddability for $\mathbb Z[\Gamma_p]$ or $\mathbb Q[\Gamma_p]$,
or for the corresponding algebras of $\Lambda_p$. Another question concerns
orderability. By Theorem~B, a left order on $G_p$ would give one on
$\Lambda_p$, answering Mal'cev's Problem~1.6 negatively in characteristic
$p$. We do not know whether such an order exists.

\subsection{A smaller orderability test}
Theorem~\ref{thm:B} reduces the orderability question to the base \(G_p\). We can go one step further. In a factorization \(M=CB\), to control the rank of \(M\) on a chosen subspace it suffices to keep track of where the \emph{right} factor \(B\) sends that subspace. This observation singles out a finitely generated subgroup of the base on which a left order would already suffice.

Put \(\beta=2(m-|\mathcal E|)+r|\mathcal E|\), the width of the factorization in (4.9), and choose \(T\subseteq\{1,\ldots,r-1\}\) of size \(t\). Use the two outer products (4.6) at each regular point and \(Z(z)I_r\) at each exceptional point. Absorb the sign in (4.9) into the left factor. Apply the normalizations (5.2), subtract row zero from the other rows on each line, and retain on that line just the rows and columns in \(T\). This gives a specified factorization over \(\F_Q[G_p]\),
\[
 M_T=C_T B_T,\qquad
 M_T=\bigoplus_i M_{i,T},\qquad M_{i,T}=(e_{i,uv}-1)_{u,v\in T},
\]
through \(\beta\) columns. Let \(P_T\le G_p\) be generated by the group elements in the supports of the entries of \(B_T\). With the choices of Section~7.1 this is a definite subgroup.

\begin{proposition}\label{prop:order-test}
If \(Nt>\beta\), the subgroup \(U_T=\langle P_T,\gamma_1,\ldots,\gamma_N\rangle\) of \(\Lambda_p\) has no faithful homomorphism to \(D^\times\) for any division ring \(D\) of characteristic \(p\). Its group algebra is a domain over every field. If \(P_T\) is left-orderable, so is \(U_T\), giving a counterexample to the field formulation of Problem~1.6 in characteristic \(p\).
\end{proposition}

\begin{proof}
Suppose \(\psi:U_T\to D^\times\) is faithful. Extend it linearly to \(\F_p[U_T]\) and induce to \(\Lambda_p\). The resulting right \(D\)-space
\[
       V=\F_p[\Lambda_p]\otimes_{\F_p[U_T]}D
\]
has one copy of \(D\) for each coset in \(\Lambda_p/U_T\). Write \(V_0\) for its identity-coset coordinate. We will compute the rank of \(M_T\) as a map \emph{from} \(V_0^{Nt}\) \emph{into} \(V^{Nt}\).

The factor \(B_T\) takes its input into \(\beta\) identity-coset coordinates, because every group element in its support lies in \(U_T\). To justify the same count with coefficients in \(\F_Q\), represent \(\F_Q\) on \(a=[\F_Q:\F_p]\) coordinates by its regular matrices, as in Lemma~4.2. The intermediate space becomes \(V_0^{a\beta}\). On the left, \(M_T\) has entries over \(\F_p\), so we get \(a\) copies of the original map. Dividing both dimensions by \(a\) gives the upper bound \(\beta\).

We now show that each line contributes \(t\). Write \(\Delta_i=\diag(e_{i,uu})_{u\in T}\) and \(C_i=\Delta_i^{-1}M_{i,T}-I_t\). By \eqref{eq:iterated-shear}, the maps \(s_i^j I_t+C_i\) on \(V_0^t\) all have the same rank as \(M_{i,T}\): \(\gamma_i\) preserves \(V_0\), and the factors on the output side are invertible on \(V^t\).

If \(s_i\notin U_T\), its coset is distinct from every \(hU_T\) with \(h\in F_i\). Indeed, \(h^{-1}s_i\in U_T\) would give \(h^{-1}\in P_T\) under the retraction, and then \(s_i\in U_T\). Every entry of \(C_i\) is supported in \(F_i\). Projecting the output of \(s_iI_t+C_i\) onto the \(s_iU_T\)-coordinate therefore gives an isomorphism on \(V_0^t\).

If \(s_i\in U_T\), project instead onto \(V_0^t\). The compressed maps are \(\psi(s_i)^j I_t+C_i^0\). The entries of \(C_i^0\) are sums of images of elements in \(F_i\cap U_T\). Thus \(\psi(s_i)\) commutes with them. Its powers are distinct by faithfulness. The eigenvalue argument of Lemma~\ref{lem:shear}, in the division subring generated by these entries and \(\psi(s_i)\), makes at least one of the \(t+1\) blocks with \(0\le j\le t\) invertible. This also gives rank \(t\).

The block diagonal map has rank \(Nt>\beta\), a contradiction. Diffuseness descends from \(\Lambda_p\) to \(U_T\), giving the domain assertion. The retraction takes \(U_T\) onto \(P_T\) and its kernel is a subgroup of the bi-orderable kernel in \eqref{eq:free-kernel}. The extension argument in the proof of Theorem~\ref{thm:B} gives the assertion on left orders.
\end{proof}

One may take \(t=\lfloor\beta/N\rfloor+1\), since (4.2) gives \(\beta<Nr/4\). This test requires an order on \(P_T\). Whether one of the eligible subgroups admits such an order remains open.

\subsection{What would suffice in characteristic zero}

For a fixed group \(H\), embedding \(\Z[H]\) in a division ring is equivalent to embedding \(\Q[H]\). In a division ring containing \(\Z[H]\), the nonzero integers are central units, so that the rational group algebra embeds by localization. There is also a simple necessary condition on a proposed integral example. If \(\F_p[H]\) embeds in a division ring \(D_p\) for infinitely many primes \(p\), take a nonprincipal ultraproduct of those division rings. Reduction of coefficients followed by the given embeddings maps \(\Z[H]\) into this ultraproduct. A nonzero finite sum has an integer coefficient that survives at all but finitely many primes, so this map is injective. Thus an obstruction for \(\Z[H]\) requires obstructions for the \emph{same} \(H\) at all but finitely many primes. This argument gives no converse. The group constructed in Section~\ref{sec:uniform} satisfies the necessary condition at every prime.

For the universal question about ordered groups one may also reduce the coefficient field within a fixed characteristic. In a companion paper we prove that if \(B[H^m]\) embeds in a division ring for every \(m\ge1\), then \(K[H]\) embeds for every extension field \(K/B\) \cite[Theorem~1.1]{EmbeddingCoeff}. Here \(H^m\) denotes the direct product of \(m\) copies of \(H\). Take \(B\) to be the prime field. Finite powers of a left-orderable group are left-orderable by lexicographic comparison. Thus in a given characteristic the assertion for \emph{all} left-orderable groups over \(B\) is equivalent to the assertion for all such groups over every field of that characteristic. The passage to powers is important: this reasoning uses embeddings for the powers of \(H\), even when the desired conclusion concerns \(H\) alone. Inversion interchanges the left and right order conventions of Problem~1.6.

Appendix~\ref{app:fixed-primes} gives another constraint specific to the matrices used here. At fixed width below full rank and with a fixed bound on the supports of the factors, a shear factorization can occur at only finitely many primes. The proof gives a computable exceptional set even if the group is allowed to vary.

\section{One central element in every positive characteristic}\label{sec:uniform}

Theorem~\ref{thm:B} gives a different group and a different word \(v_p\) for each prime \(p\). Their restricted direct product contains all the words. In characteristic \(p\), the theorem forces the word in the \(p\)-factor to die. We put this product inside a larger group and arrange that the commutators of all the \(v_p\) with one new generator have the same value.

\begin{maintheorem}\label{thm:D}
There is a three-generated diffuse group \(E\) with a central element \(z\) of infinite order such that the following hold.
\begin{enumerate}
\item Every homomorphism \(E\to D^\times\) sends \(z\) to \(1\) whenever \(D\) is a division ring of positive characteristic.
\item Every finite-dimensional linear representation of \(E\) over any field sends \(z\) to the identity. More generally, every homomorphism from \(E\) to a group satisfying the descending chain condition on centralizers of finite subsets kills \(z\).
\end{enumerate}
\end{maintheorem}

\subsection{Keeping the number of generators finite}

We first embed the direct product with finite support in a group with just two generators. The embedding uses two unrestricted wreath products, a device going back to B.~H.~Neumann and H.~Neumann \cite{NeumannNeumann}; the first embedding of countable groups in two-generated groups, by Higman, B.~H.~Neumann and H.~Neumann \cite{HNN}, used HNN extensions instead. We follow the version of Darbinyan \cite[Section~3]{Darbinyan}. We include the construction to check diffuseness and to keep track of the embedded subgroups.

\begin{lemma}\label{lem:uniform-packing}
Every countable diffuse group embeds in a two-generated diffuse group.
\end{lemma}
\begin{proof}
Finite direct products of diffuse groups are diffuse by Lemma~\ref{lem:extension}. For an arbitrary product, take a finite subset \(S\) with at least two elements. Some finite list of coordinates distinguishes its elements. Project onto these coordinates and lift two extremal points of the image to \(S\). If \(x\) is either lift and both \(hx,h^{-1}x\) belong to \(S\), extremality forces \(h\) to project to the identity. Then \(hx\) and \(x\) have the same projection, so \(hx=x\). This proves diffuseness of the product. Consequently \(C^{\Z}\rtimes\Z\) is diffuse whenever \(C\) is diffuse, again by Lemma~\ref{lem:extension}.

Let \(a_1,a_2,\ldots\) generate the countable diffuse group \(A\). In \(A^{\Z}\rtimes\langle t\rangle\) use the shift \((tft^{-1})(n)=f(n-1)\). Let \(f_i(n)=a_i\) for \(n\ge0\) and \(f_i(n)=1\) for \(n<0\). Then \([f_i,t]=\delta_0(a_i)\), where \(\delta_j(a)\) has value \(a\) at \(j\) and the identity elsewhere. Put \(H=\langle t,f_1,f_2,\ldots\rangle\).

In \(H^\Z\rtimes\langle u\rangle\) use the same shift and define a function \(c\) by \(c(1)=t\), \(c(2^i)=f_i\) for \(i\ge1\), and \(c(n)=1\) elsewhere. For \(d_i=2^i-1\), the supports of \(u^{-d_i}cu^{d_i}\) and \(c\) can meet only at \(1\). Indeed an overlap requires \(2^a-2^b=2^i-1\) for some \(a>b\ge0\). The right side is odd, so \(b=0\) and \(a=i\). At the common coordinate the two values are \(f_i\) and \(t\). Thus
\[
 [u^{-d_i}cu^{d_i},c]=\delta_1([f_i,t])=\delta_1(\delta_0(a_i)).
\]
These elements generate a coordinate copy of \(A\) in \(B=\langle u,c\rangle\). Both \(H\) and \(B\) lie in diffuse wreath products, and diffuseness passes to subgroups.
\end{proof}

\subsection{Making all the commutators equal}

\begin{proof}[Proof of Theorem~\ref{thm:D}]
Take \(A=\bigoplus_p\Lambda_p\), the subgroup of the direct product consisting of tuples with finite support. It is countable and diffuse by Theorem~\ref{thm:B} and the product argument above, since the restricted product is a subgroup of the unrestricted one. Embed \(A\) in a two-generated diffuse group \(B\) by Lemma~\ref{lem:uniform-packing}. Each \(v_p\) from Theorem~\ref{thm:B} has infinite order because \(\Lambda_p\) is diffuse. These elements belong to different direct factors, so \(J=\bigoplus_p\langle v_p\rangle\) is free abelian in \(B\). Let \(\chi:J\to\Z\) send each \(v_p\) to \(1\).

Form the abelian group \(M\) with generators \(z\) and \(x_b\) for \(b\in B\), and relations
\begin{equation}\label{eq:uniform-module}
 x_{bj}=x_b+\chi(j)z \qquad(b\in B,\ j\in J).
\end{equation}
Here \(M\) is written additively. Choose one representative \(t\) from each left coset \(bJ\), and let \(T\) be this set of representatives. Then
\[
 M=\Z z\ \oplus\!\bigoplus_{t\in T}\Z x_t .
\]
Indeed each \(b\) has a unique expression \(tj\). Sending \(x_{tj}\) to \(x_t+\chi(j)z\) and \(z\) to itself respects \eqref{eq:uniform-module} since \(\chi\) is a homomorphism. This map inverts the natural map from the displayed free abelian group to \(M\). In particular \(z\) has infinite order. Left translation \(x_b\mapsto x_{gb}\) for \(g\in B\) respects the same relations and defines an action of \(B\) fixing \(z\).

Put \(E=M\rtimes B\). The element \(z\) is central, and conjugating \(x_1\) by \(b\) gives \(x_b\). For every prime \(p\) the commutator in \(E\), evaluated in \(M\), is
\begin{equation}\label{eq:uniform-commutator}
 [v_p,x_1]=x_{v_p}-x_1=z .
\end{equation}
Thus the two generators of \(B\) and \(x_1\) generate \(E\). The free abelian group \(M\) is diffuse, so Lemma~\ref{lem:extension} gives diffuseness of \(E\). A homomorphism \(E\to D^\times\) with \(\operatorname{char}D=p\) kills \(v_p\) on restriction to \(\Lambda_p\), by Theorem~\ref{thm:B}. Equation~\eqref{eq:uniform-commutator} then kills \(z\).

Now suppose a homomorphism \(\rho:E\to H\) satisfies \(\rho(z)\ne1\). Equation~\eqref{eq:uniform-commutator} forces every \(\rho(v_p)\) to be nontrivial. Since \(v_p\in[\Lambda_p,\Lambda_p]\), each subgroup \(\rho(\Lambda_p)\) is nonabelian. The images for different primes commute elementwise because their factors commute in \(A\). Enumerate the primes \(p_1,p_2,\ldots\) and choose noncommuting \(a_i,b_i\in\rho(\Lambda_{p_i})\). Then \(b_i\) centralizes \(a_1,\ldots,a_{i-1}\) and fails to centralize \(a_i\). Hence
\[
 H\supsetneq C_H(a_1)\supsetneq C_H(a_1,a_2)\supsetneq\cdots .
\]
This contradicts the stated chain condition. For \(H=\mathrm{GL}_n(K)\), the matrices in \(M_n(K)\) that commute with a given finite set form a vector subspace, whose intersection with \(H\) is the group centralizer. Each strict inclusion above forces a strict drop in the dimension of these subspaces. There can be at most \(n^2\) such drops, over any field \(K\).
\end{proof}

\begin{corollary}\label{cor:uniform-localization}
For every field \(K\) of positive characteristic, \(z-1\) is a nonzero central element of the domain \(K[E]\) killed by every unital homomorphism to a division ring. The central localization \(K[E][(z-1)^{-1}]\) is a nonzero domain admitting no unital homomorphism to a division ring.
\end{corollary}
\begin{proof}
Diffuseness gives the domain assertion, and \(z\) and \(1\) are distinct group basis elements. A homomorphism \(K[E]\to D\) to a division ring restricts to a group homomorphism \(E\to D^\times\), and \(D\) has the characteristic of \(K\). Theorem~\ref{thm:D} sends \(z-1\) to zero. Its powers are nonzero central elements, so localization preserves the domain. A map from the localization would have to send \(z-1\) to a unit.
\end{proof}

Domains with no homomorphism to a division ring have been known since Leavitt and Cohn, for instance Leavitt's algebras \(V_{m,n}\) with \(2\le m<n\) \cite[Theorem~6.1 and Example~6.3]{Bergman}; the point here is that this one is a central localization of a group algebra.

The same \(z\) dies under every homomorphism to a bi-orderable group \(P\), since \(P\) embeds by monomials in the Mal'cev--Neumann division ring \(\F_2((P))\) \cite{Neumann}\cite[Section~1]{r1}. In particular \(E\) admits no bi-order. Whether a division ring of characteristic zero can detect \(z\) remains open.

\appendix
\section{Formal statements and their scope}

The Lean~4 development in the archived release \cite{r18} proves the
nonembedding clauses of Theorems~A(ii) and~B. It proves existence of
admissible data for each prime. For every such datum $d$, every field $K$
of characteristic $p$, and every division ring $D$, no unital homomorphism
into $D$ from either $K[d.\mathrm{Gamma}]$ or $K[d.\mathrm{Shear}]$ is
injective. The following table uses the common prefix
\texttt{DomainNonembedding.} and abbreviates its namespaces
\texttt{Main}, \texttt{Permutation} and \texttt{ShearNonembedding} by
\texttt{M}, \texttt{P} and \texttt{S}.

\begin{center}
\begingroup
\tabsmall
\setlength{\tabcolsep}{0pt}
\begin{tabular}{@{}p{.425\linewidth}@{\hspace{.015\linewidth}}p{.55\linewidth}@{}}
\textbf{Lean declaration} & \textbf{Mathematical statement}\\
\hline
\texttt{M.paperData\_exists}
  & Lemma~4.1 and the two-block partition.\\
\texttt{M.interpolation\_identity}
  & The punctured-line identity (4.7).\\
\texttt{M.rank\_certificate}
  & Lemma~4.2.\\
\texttt{M.universal\_rank\_ceiling}
  & Proposition~4.3, including strictness.\\
\texttt{M.line\_free}
  & The free line basis in Lemma~5.2.\\
\texttt{M.central\_injective},\newline
\texttt{M.copied\_injective}
  & The vertex injections in Section~3.\\
\texttt{P.average},\newline
\texttt{P.exists\_common\_normalized}
  & Lemma~3.1 and one permutation for all lines.\\
\texttt{M.nonembedding\_prime\_data}
  & Nonembedding for every datum over the prime field.\\
\texttt{M.nonembedding\_data}
  & Theorem~A(ii), every admissible datum and field of characteristic $p$.\\
\texttt{M.nonembedding}
  & The same conclusion for a selected witness for each prime.\\
\texttt{S.nonembedding\_data}
  & Theorem~B, nonembedding for every datum and field of characteristic $p$.\\
\texttt{M.shear\_nonembedding}
  & The shear conclusion for a selected witness.\\
\texttt{S.shear\_rank\_bound}
  & Half-rank shear bound sufficient for Theorem~B; not Lemma~2.2 in full.
\end{tabular}
\endgroup
\end{center}

\noindent\textit{Comparing the presentations.}
The structure \texttt{PaperData} and the predicate \texttt{Admissible}
specify the field and its numbering, distinct marks, nonzero directions,
punctured incidences, girth, density and exceptional-set bounds, and two
blocks with bijective coordinates at each regular point. The deterministic
datum of Section~7 satisfies these conditions. The all-data theorems apply
to it directly.

For this datum, \texttt{d.BaseGroup} uses the triangles (4.4) and line
words (5.2). Its presentation collapses row and column zero of one
distinguished line and one edge $L_{b,0,0}$ at each apex. The tree chosen
in Section~7 may differ. Both presentations compute the fundamental group
of the same based complex. Changing the base paths transports each line
basis by a common conjugation. Equation (3.3) uses one stable letter for an
entire edge basis and kills precisely the spanning-tree letters $t_{\pi,i_*}$.
Hence \texttt{d.Gamma} presents the permutation group, and \texttt{d.Shear}
uses precisely (7.1) on the same base. The cones in Section~6 impose no
relations on either group. These identifications are a comparison on paper.
The separate witnesses selected by choice in \texttt{paperData} are covered
by the formal theorems for selected witnesses, without a claim that their
numbered presentations coincide with ours.

\medskip
\noindent\textit{Rank and geometry.}
The imported primitive \texttt{NonSofic.rankD} uses left row rank. Its
equality with right column rank (transpose over the opposite division ring)
is proved by the two lemmas
\begin{itemize}[itemsep=0pt,topsep=2pt]
\item \texttt{DomainNonembedding.\allowbreak DivisionRingRank.\allowbreak row\_eq\_column},
\item \texttt{DomainNonembedding.\allowbreak rowRank\_eq\_columnRank}.
\end{itemize}
Both the primitive and the bridge are included in the axiom audit. The
formal shear proof applies the argument of (2.4) to the matrix $(\rho(e_{uv}))$ of group entries. A map to an auxiliary HNN construction proves the
generators $s_i$ are nonidentity. The formal argument requires only that
fact about their images. The full-rank Lemma~2.2 and the general attachment
Theorem~C have written proofs. The geometric and domain assertions, finite
presentability, torsion-freeness, the retraction and orderability statements
likewise rely on the written proofs. The formalization makes no
characteristic-zero assertion.

The common words in Theorems~A(iii) and~B, the orderability test in Proposition~\ref{prop:order-test}, the central element of Theorem~D, and the finite-prime bound of Appendix~\ref{app:fixed-primes} have written proofs. The table records the scope of the archived Lean development.

\medskip
\noindent\textit{Reproduction.}
Release \texttt{v1.1.0} pins Lean~4.35.0 and Mathlib commit
\texttt{56a005855ea17e0b0958f7\allowbreak 3ec801ae140b0a1679}.
From its \texttt{lean/} directory, run \texttt{lake exe cache get}, then
\texttt{lake build}. Its \texttt{README.md}, \texttt{SCOPE.md},
\texttt{REPORT.md} and \texttt{evidence/} record the audit and replay
commands. The archived default build comprised 3657 jobs. The audit
covered 2792 declarations and three release wrappers, with only
\texttt{propext}, \texttt{Classical.choice} and \texttt{Quot.sound}.
There are no admissions in the development or its default target. A
separate challenge contains three intentional template holes. The archived
checks include a sandboxed comparator against frozen statements and
successful replay of the sixteen-root proof closure with \texttt{leanchecker}
and \texttt{nanoda}. Removing \texttt{Classical.choice} from the allowed
axioms failed as expected. The clean release reproduced the recorded audit
and proof-closure hash.

\section{Finitely many primes at fixed size}\label{app:fixed-primes}
The full-rank calculation for a shear also holds in the ring of operators affiliated to the complex group von Neumann algebra. It therefore limits how far a fixed-size group-ring factorization can carry the positive-characteristic proof.

\begin{proposition}\label{prop:bounded-shear-primes}
Fix sizes \(n_1,\ldots,n_t\geq 1\), \(a\geq 0\), width \(0\leq b<a+\sum_i n_i\), and support bound \(h\geq 1\). Let \(K\) be a field and \(H\) any group containing arrays \((e_{i,uv})_{u,v=1}^{n_i}\) and elements \(s_i,\gamma_i\) satisfying
\[
 [s_i,\gamma_i]=[s_i,e_{i,uv}]=1,\qquad
 \gamma_i e_{i,uv}\gamma_i^{-1}=s_i^{\delta_{uv}}e_{i,uv},
 \qquad \operatorname{ord}(s_i)=\infty .
\]
Put \(M_i=(e_{i,uv}-1)_{u,v}\) and \(T=\operatorname{diag}(I_a,M_1,\ldots,M_t)\). Suppose \(T=PQ\) over \(K[H]\), where \(P\) has \(b\) columns, \(Q\) has \(b\) rows, and every entry of either factor has at most \(h\) support elements. This is impossible in characteristic zero. In characteristic \(p>0\), it can occur only if \(p\) divides a computable positive integer depending just on \(a,n_1,\ldots,n_t,b,h\). Thus a single finite set of primes works even when the group, its elements, and the coefficient field vary.
\end{proposition}
We use von Neumann dimension for the regular action of a discrete group \(H\) on \(\ell^2(H)\), with \(\dim_H\ell^2(H)=1\). In the ring of affiliated operators every finite matrix has a partial inverse. Its image and kernel are therefore finitely generated projectives, on which dimension is additive and faithful. For a subgroup the ranks of its matrices are preserved on passage to \(H\), by decomposing the regular representation into cosets \cite[Chapters~6 and~8]{LuckAff}.

\begin{proof}
For one block \(M_i\), work in the ring of affiliated operators of the subgroup \(L_i=\langle s_i,\gamma_i,e_{i,uv}\rangle\). Write \(s=s_i\), \(M=M_i\), \(\Delta=\operatorname{diag}(e_{i,uu})_u\), and \(C=\Delta^{-1}(M-\Delta)\). As in the shear calculation, conjugating \(M\) by \(\gamma_i^j\) gives \(\Delta(s^jI+C)\). Hence the kernels of \(s^jI+C\) all have equal von Neumann dimension. They form an algebraic direct sum. Indeed \(s\) is central in this affiliated ring, and applying \(\prod_{\ell\ne j}(s^\ell I+C)\) to a relation among finitely many of these kernels kills all but the \(j\)-th term. On that term it acts by a product of \(s^\ell-s^j\), each invertible because a nontrivial power of \(s\) has infinite orbits on \(L_i\) and fixes no nonzero vector in \(\ell^2(L_i)\). Arbitrarily many of the kernels fit in a module of dimension \(n_i\), so each has dimension zero. Subgroup inclusion preserves affiliated rank, and \(T\) thus has full rank \(a+\sum_i n_i\) over the complex affiliated ring of \(H\). The factors would bound its rank by \(b\). Any characteristic-zero factorization could be read over the complex numbers by embedding the finitely generated field of its coefficients in \(\mathbb C\).

To obtain the uniform set of primes, pad the entries of \(P,Q\) to \(h\) formal slots with coefficient variables. In each entry of \(PQ-T\), record the partition of the finitely many formal terms by equality of their group words. Each part gives an integer polynomial equation by summing its coefficients. There are only finitely many partition patterns at the prescribed dimensions and support bounds. For any pattern realized by group words satisfying the shear hypotheses, the polynomial equations have no complex solution, since such a solution would give the same factorization in that group over \(\mathbb C\). The Nullstellensatz makes the ideal the unit ideal over \(\mathbb Q\). For every enumerated pattern with unit ideal over \(\mathbb Q\), a Gr\"obner basis computation with certificates expresses \(1\) in the ideal. Clearing denominators gives a positive integer in the corresponding ideal over \(\mathbb Z\). Take the product of these integers, or \(1\) if the list is empty. A solution in characteristic \(p\) makes the integer for its pattern vanish, so \(p\) divides the product.
\end{proof}

The sizes and supports in this statement are those of an actual factorization over the group algebra, after any matrix amplification used to represent extension-field scalars. The bound places no restriction on proofs in which these sizes grow with the prime.

\section*{Acknowledgements}

We thank Andrei Jaikin-Zapirain for his questions about embeddings of
group algebras and for the correspondence that prompted this work.

Claude Opus 5.5, GPT-6.1 Sol, GPT-6.2 Sol and GPT-6.2 Astra assisted with this work; the authors take responsibility for every claim. The Lean development and its verification records are archived with the paper~\cite{r18}.

\end{document}